\documentclass[12pt]{article}

\def\ontop#1#2{\buildrel \rm #1 \over #2}

\def\qed{\hfill\ \rule{2mm}{2mm} }

\def\mat{\bf }

\newcounter{itemtransfer}
\def\bibart#1#2#3#4#5#6#7
{\bibitem{#1} #2, #4, #5 #6 (#3), #7.}

\def\bibbook#1#2#3#4#5#6
{\bibitem{#1}
#2, #4, (#5, #6, #3).}

\def\bibcoll#1#2#3#4#5#6#7#8
{\bibitem{#1}
#2, #4, in: #5, #6 (#7, #3) p. #8.}

\def\bibdiss#1#2#3#4#5#6
{\bibitem{#1}
#2, #4, #5, #3.}

\newtheorem{guess}{Guess}[section]

\newtheorem{define}[guess]{Definition}
\newtheorem{prop}[guess]{Proposition}
\newtheorem{theorem}[guess]{Theorem}

\newtheorem{lem}[guess]{Lemma}

\newtheorem{remark}[guess]{Remark}

\begin{document}

\bibliographystyle{plain}

\title{
Nonisomorphic Ordered Sets with Arbitrarily Many Ranks That
Produce Equal Decks
\thanks{This work was sponsored in part by
Louisiana Board of Regents RCS grant LEQSF(1999-02)-RD-A-27.
}
}

\author{ \small Bernd S. W. Schr\"oder \\
\small Program of Mathematics \& Statistics\\
\small Louisiana Tech University\\
\small Ruston, LA 71272\\
\small {\tt schroder@coes.LaTech.edu }
}

\date{\small \today}

\maketitle

\begin{abstract}

\normalsize
We prove that for any $n$ there is a pair $(P_1 ^n , P_2 ^n )$ of nonisomorphic ordered sets
such that $P_1 ^n $ and $P_2 ^n $
have equal maximal and minimal decks, equal neighborhood decks, and there are
$n+1$ ranks $k_0 , \ldots , k_n $
such that for each $i$ the decks obtained by removing the
points of rank $k_i $ are equal.
The ranks $k_1 , \ldots , k_n $ do not contain extremal elements and at each of the other ranks there are
elements whose removal will produce isomorphic cards.
Moreover, we show that such sets can be constructed such that
only for ranks $1$ and $2$, both without extremal elements,
the decks obtained by removing the
points of rank $r_i $ are not equal.

\end{abstract}


\section{Introduction}

The reconstruction problem for ordered sets asks if it is possible to reconstruct the isomorphism type of
a given ordered set from its collection
(the ``deck")
of one-point-deleted subsets.
In \cite{Sanun}, Sands asked
if ordered sets might even be reconstructible from
the the collection of subsets obtained by erasing single maximal elements (the ``maximal deck").
The negative answer to Sands' question in
\cite{KRex}
together with the paper \cite{KRtow} were the starting point for serious investigation of
the reconstruction problem for ordered sets.
(Reconstruction of graphs and other relations has a
longer history, see \cite{BoHem,Manv,Pouinterdit,Stockwindmill}.)
Since then, results that show reconstructibility given
certain types of information (for example, see \cite{IlRaCG}) as well as
results on reconstructibility of certain classes of ordered sets (for example, see \cite{RaSch})
have been proved. For a more comprehensive survey of results
and references available to date, consider \cite{JXsurvey,Schbook}.

Recently in \cite{Schexam} it has been shown that even the maximal and minimal decks together
are not sufficient to reconstruct ordered sets.
Moreover, it was shown in \cite{Schexam} that there are families of $\displaystyle{ 2^{O(\sqrt[3]{n} )} } $
pairwise nonisomorphic ordered sets of size $n$ that
all have equal maximal and minimal decks.

The desire to
rely on only a limited, focused amount of information in a reconstruction proof,
and the fact that almost all ordered sets are reconstructible from two identifiable maximal cards
(see \cite{Schrigid}, Corollary 3.10)
motivates extensions of Sands' question.
What subsets of the deck have a reasonable chance to
effect reconstruction? In \cite{Schmoreexam} it was shown that the maximal deck plus the minimal deck plus one deck
obtained by removing points of a rank $k$ that contains no extremal elements are not sufficient for
reconstruction.
The examples in \cite{Schmoreexam}
are somewhat limited. They
could not be extended to an ambiguity with more than two sets.
It also appeared as if they could
not be extended to more than two maximal elements, two minimal elements and two elements of rank $k$
or to more than one middle rank producing equal decks.
Immediately two questions arise.
\begin{itemize}
\item
Are ordered sets reconstructible from the maximal deck, the minimal deck and a deck obtained
by removing points of rank $k$ if one of these decks has at least three cards?

\item
Are ordered sets reconstructible from the maximal deck, the minimal deck and
two decks that were obtained by removing points of ranks $k$ and $l$ with $k\not= l $?

\end{itemize}

In this paper these two questions are answered negatively, even if the neighborhood
decks are equal and for any two isomorphic cards the neighborhoods of the
removed elements are also isomorphic.
The present examples provide new guidance as to what kind of
partial information is at least needed to
reconstruct ordered sets.
In particular, they show that information derived from
``small" ranks, and even from many small ranks, is not sufficient
to effect reconstruction.
Analysis of the examples also leads to results that underscore the role of rigidity in
order reconstruction (cf. Section \ref{rigsect}).
Ideas on what types of information to consider next are given in the conclusion.

\section{Basic Definitions and Preliminaries}

An {\bf ordered set} is a set $P$ equipped with a reflexive, antisymmetric and
transitive relation $\leq $, the order relation.
Throughout this paper we will assume that all ordered sets involved are finite.
Elements $x,y\in P$ are called {\bf comparable} iff $x\leq y$ or $y\leq x$.
An {\bf antichain} is an ordered set in which each element is only comparable to itself.
A {\bf chain} is an ordered set
in which any two elements are comparable.
The {\bf length} of a chain is its number of elements minus $1$.
An element $m\in P$ is called {\bf maximal} iff for all $x$ comparable to $m$ we
have $x\leq m$. {\bf Minimal} elements are defined dually.
The {\bf rank} of an element $x\in P$ is the length of the longest chain that has a minimal element
as its smallest element and $x$ as its largest element.

The {\bf dual} $P^d $ of an ordered set $P$ is the ordered set obtained by reversing all
comparabilities.
The {\bf dual rank} of an element $x\in P$ is the length of the longest chain that has a maximal element
as its largest element and $x$ as its smallest element.

A function $f:P\to Q$ from the ordered set $P$ to the ordered set $Q$ is called {\bf order-preserving}
iff for all $x,y\in P$ we have that $x\leq y$ implies $f(x)\leq f(y)$.
The function $\varphi :P\to Q$ is called an {\bf (order) isomorphism} iff $\varphi $ is bijective, order-preserving
and $\varphi ^{-1} $ is order-preserving, too.
An order isomorphism with equal domain and range is called an {\bf (order) automorphism}.
An ordered set with exactly one order automorphism (the identity)
is called {\bf rigid}.


For precise overall reconstruction terminology, cf. \cite{Schbook,Schrigid,Schmoreexam}.
For the purposes of this paper, a {\bf card} of an ordered set is a subset
with one point deleted.
If the deleted element is of rank $k$ we shall also call the card a {\bf rank $k$ card}.
The set of all cards obtained by erasing elements of rank $k$ is called the
{\bf rank $k$ deck}.
The set of all rank $k$ decks is called the {\bf ranked deck}.
A rank $k$ card is {\bf marked} iff there is a function that indicates the
rank of each element in the original set.
The set of all marked rank $k$ cards is the {\bf marked rank $k$ deck} and the
set of all marked rank $k$ decks is the {\bf marked ranked deck}.
A {\bf maximal card} is a subset in which a maximal element is erased
and a {\bf minimal card} is a subset in which a minimal element is erased.
The sets of maximal and minimal cards respectively are called the {\bf maximal deck} and the
{\bf minimal deck}.
{\bf Marked maximal cards} are maximal cards for which there is a function that indicates which
elements are maximal in the original set. The set of all marked maximal cards is called the
{\bf marked maximal deck}.
{\bf Marked minimal cards} and the {\bf marked minimal deck} are defined dually.
Isomorphic cards will also be called {\bf equal cards}, because their isomorphism
classes are equal. Decks will be called {\bf equal} iff there is a
bijection such that each card is isomorphic to its image.
Marked cards will be called equal iff there is an isomorphism that preserves the marked property
(rank in the original set or maximality/minimality in the original set). Marked decks will be called
equal iff there is a bijection such that any set is isomorphic to its image and each
isomorphism also preserves the marked property.
The {\bf up-set} of an element $x$ is the set $\uparrow x=\{ p\in P:p\geq x\} $
and the {\bf down-set} is $\downarrow x=\{ p\in P:p\leq x\} $.
The {\bf neighborhood} of an element $x$ is the set
$\updownarrow x=\uparrow x\cup \downarrow x$.
The set of all neighborhoods of points of rank $k$ is called the
{\bf rank $k$ neighborhood deck}.

We are concerned with results that show what type of information is {\em not}
sufficient to effect reconstruction. Therefore, throughout we
construct nonisomorphic ordered sets such that between some of their cards
there are isomorphisms with certain
properties.

\section{Pairs of Nonisomorphic Ordered Sets with Equal Maximal and Minimal Decks and
$n+1$ Non-Extremal Ranks for Which the Rank $k$ Decks are Equal}

In this section we describe the fundamental construction used to build the examples.
Lemma \ref{severalmiddle} gives the overall idea, which is an extension of the work in \cite{Schmoreexam}.
Lemma \ref{severalmiddle} is also quite similar to
the examples in \cite{Pouinterdit}.
In a way, Lemma \ref{severalmiddle} is reminiscent of a ``vertical M\"obius strip".

Lemma \ref{getsmallerR} then shows that sets as needed in the construction in
Lemma \ref{severalmiddle} actually exist.
In the following, we explore some features of the construction as well as variations
that lead to examples with
other properties.

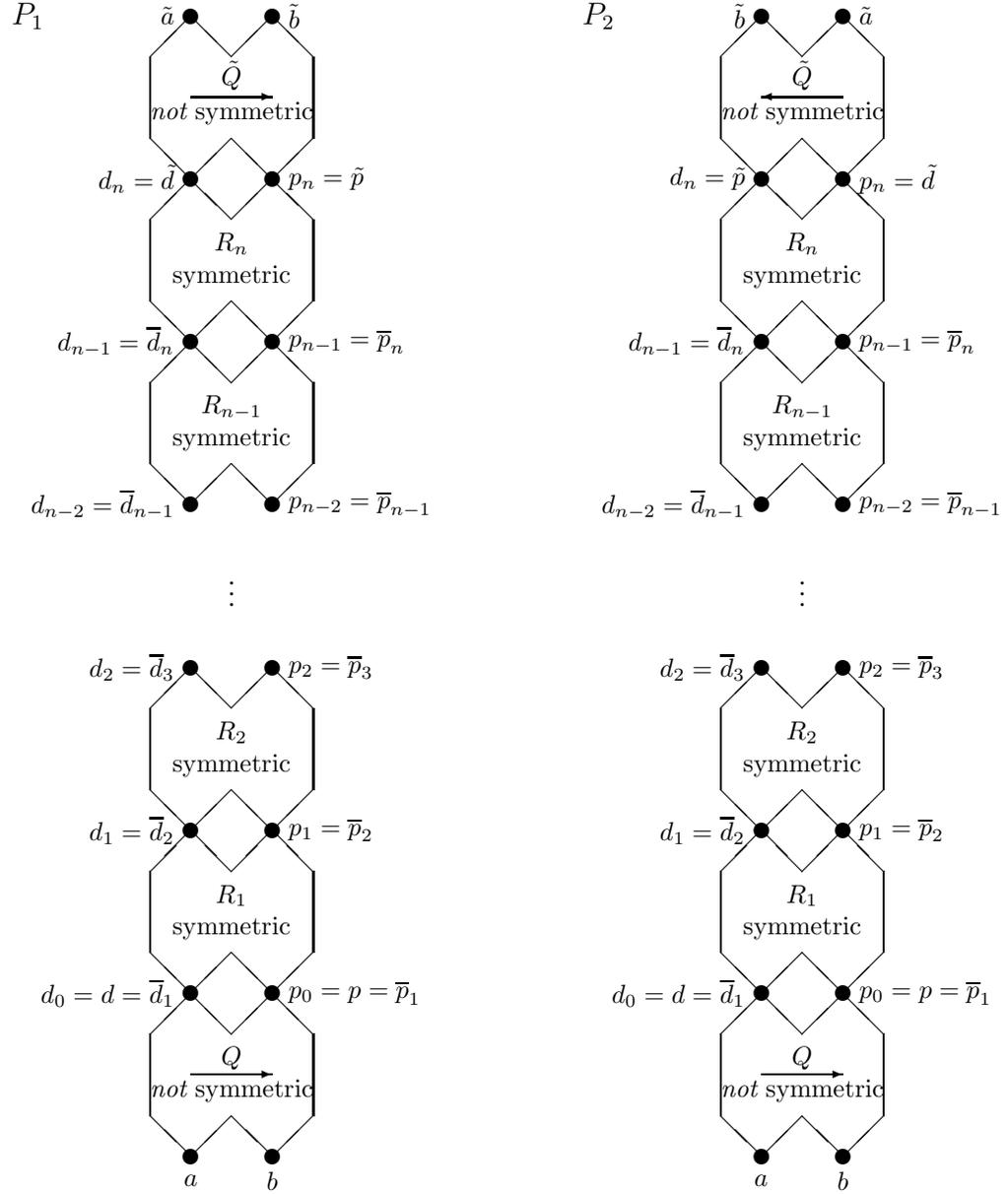
\begin{figure}

\centerline{
\unitlength 1.1mm
\linethickness{0.4pt}
\begin{picture}(110.00,146.00)
\put(20.00,10.00){\line(1,-1){5.00}}
\put(25.00,5.00){\line(1,1){5.00}}
\put(30.00,10.00){\line(1,-1){5.00}}
\put(35.00,5.00){\line(1,1){5.00}}
\put(40.00,10.00){\line(0,1){10.00}}
\put(40.00,20.00){\line(-1,1){5.00}}
\put(35.00,25.00){\line(-1,-1){5.00}}
\put(30.00,20.00){\line(-1,1){5.00}}
\put(25.00,25.00){\line(-1,-1){5.00}}
\put(20.00,20.00){\line(0,-1){10.00}}
\put(30.00,17.00){\makebox(0,0)[cc]{\footnotesize $Q$}}
\put(30.00,13.00){\makebox(0,0)[cc]{\footnotesize {\em not} symmetric}}
\put(25.00,15.00){\vector(1,0){10.00}}
\put(35.00,5.00){\circle*{2.00}}
\put(25.00,5.00){\circle*{2.00}}
\put(25.00,2.00){\makebox(0,0)[cc]{\footnotesize $a$}}
\put(35.00,2.00){\makebox(0,0)[cc]{\footnotesize $b$}}
\put(37.00,25.00){\makebox(0,0)[lc]{\footnotesize $p_0 = p=\overline{p}_1 $}}
\put(23.00,25.00){\makebox(0,0)[rc]{\footnotesize $d_0 = d=\overline{d}_1 $}}
\put(20.00,30.00){\line(1,-1){5.00}}
\put(25.00,25.00){\line(1,1){5.00}}
\put(30.00,30.00){\line(1,-1){5.00}}
\put(35.00,25.00){\line(1,1){5.00}}
\put(40.00,30.00){\line(0,1){10.00}}
\put(40.00,40.00){\line(-1,1){5.00}}
\put(35.00,45.00){\line(-1,-1){5.00}}
\put(30.00,40.00){\line(-1,1){5.00}}
\put(25.00,45.00){\line(-1,-1){5.00}}
\put(20.00,40.00){\line(0,-1){10.00}}
\put(30.00,37.00){\makebox(0,0)[cc]{\footnotesize $R_1 $}}
\put(30.00,33.00){\makebox(0,0)[cc]{\footnotesize symmetric}}
\put(35.00,25.00){\circle*{2.00}}
\put(25.00,25.00){\circle*{2.00}}
\put(37.00,45.00){\makebox(0,0)[lc]{\footnotesize $p_1 =\overline{p}_2 $}}
\put(23.00,45.00){\makebox(0,0)[rc]{\footnotesize $d_1=\overline{d}_2 $}}
\put(20.00,50.00){\line(1,-1){5.00}}
\put(25.00,45.00){\line(1,1){5.00}}
\put(30.00,50.00){\line(1,-1){5.00}}
\put(35.00,45.00){\line(1,1){5.00}}
\put(40.00,50.00){\line(0,1){10.00}}
\put(40.00,60.00){\line(-1,1){5.00}}
\put(35.00,65.00){\line(-1,-1){5.00}}
\put(30.00,60.00){\line(-1,1){5.00}}
\put(25.00,65.00){\line(-1,-1){5.00}}
\put(20.00,60.00){\line(0,-1){10.00}}
\put(30.00,57.00){\makebox(0,0)[cc]{\footnotesize $R_2 $}}
\put(30.00,53.00){\makebox(0,0)[cc]{\footnotesize symmetric}}
\put(35.00,45.00){\circle*{2.00}}
\put(25.00,45.00){\circle*{2.00}}
\put(37.00,65.00){\makebox(0,0)[lc]{\footnotesize $p_2 =\overline{p}_3 $}}
\put(23.00,65.00){\makebox(0,0)[rc]{\footnotesize $d_2=\overline{d}_3 $}}
\put(35.00,65.00){\circle*{2.00}}
\put(25.00,65.00){\circle*{2.00}}
\put(30.00,75.00){\makebox(0,0)[cc]{$\vdots $}}
\put(37.00,85.00){\makebox(0,0)[lc]{\footnotesize $p_{n-2} =\overline{p}_{n-1} $}}
\put(23.00,85.00){\makebox(0,0)[rc]{\footnotesize $d_{n-2} =\overline{d}_{n-1} $}}
\put(20.00,90.00){\line(1,-1){5.00}}
\put(25.00,85.00){\line(1,1){5.00}}
\put(30.00,90.00){\line(1,-1){5.00}}
\put(35.00,85.00){\line(1,1){5.00}}
\put(40.00,90.00){\line(0,1){10.00}}
\put(40.00,100.00){\line(-1,1){5.00}}
\put(35.00,105.00){\line(-1,-1){5.00}}
\put(30.00,100.00){\line(-1,1){5.00}}
\put(25.00,105.00){\line(-1,-1){5.00}}
\put(20.00,100.00){\line(0,-1){10.00}}
\put(30.00,97.00){\makebox(0,0)[cc]{\footnotesize $R_{n-1} $}}
\put(30.00,93.00){\makebox(0,0)[cc]{\footnotesize symmetric}}
\put(35.00,85.00){\circle*{2.00}}
\put(25.00,85.00){\circle*{2.00}}
\put(37.00,105.00){\makebox(0,0)[lc]{\footnotesize $p_{n-1} =\overline{p}_n $}}
\put(23.00,105.00){\makebox(0,0)[rc]{\footnotesize $d_{n-1} =\overline{d}_n $}}
\put(20.00,110.00){\line(1,-1){5.00}}
\put(25.00,105.00){\line(1,1){5.00}}
\put(30.00,110.00){\line(1,-1){5.00}}
\put(35.00,105.00){\line(1,1){5.00}}
\put(40.00,110.00){\line(0,1){10.00}}
\put(40.00,120.00){\line(-1,1){5.00}}
\put(35.00,125.00){\line(-1,-1){5.00}}
\put(30.00,120.00){\line(-1,1){5.00}}
\put(25.00,125.00){\line(-1,-1){5.00}}
\put(20.00,120.00){\line(0,-1){10.00}}
\put(30.00,117.00){\makebox(0,0)[cc]{\footnotesize $R_n $}}
\put(30.00,113.00){\makebox(0,0)[cc]{\footnotesize symmetric}}
\put(35.00,105.00){\circle*{2.00}}
\put(25.00,105.00){\circle*{2.00}}
\put(37.00,125.00){\makebox(0,0)[lc]{\footnotesize $p_n =\tilde{p} $}}
\put(23.00,125.00){\makebox(0,0)[rc]{\footnotesize $d_n=\tilde{d} $}}
\put(35.00,125.00){\circle*{2.00}}
\put(25.00,125.00){\circle*{2.00}}
\put(20.00,130.00){\line(1,-1){5.00}}
\put(25.00,125.00){\line(1,1){5.00}}
\put(30.00,130.00){\line(1,-1){5.00}}
\put(35.00,125.00){\line(1,1){5.00}}
\put(40.00,130.00){\line(0,1){10.00}}
\put(40.00,140.00){\line(-1,1){5.00}}
\put(35.00,145.00){\line(-1,-1){5.00}}
\put(30.00,140.00){\line(-1,1){5.00}}
\put(25.00,145.00){\line(-1,-1){5.00}}
\put(20.00,140.00){\line(0,-1){10.00}}
\put(30.00,137.67){\makebox(0,0)[cc]{\footnotesize $\tilde{Q} $}}
\put(30.00,133.00){\makebox(0,0)[cc]{\footnotesize {\em not} symmetric}}
\put(25.00,135.00){\vector(1,0){10.00}}
\put(37.00,145.00){\makebox(0,0)[lc]{\footnotesize $\tilde{b} $}}
\put(23.00,145.00){\makebox(0,0)[rc]{\footnotesize $\tilde{a} $}}
\put(35.00,145.00){\circle*{2.00}}
\put(25.00,145.00){\circle*{2.00}}
\put(5.00,145.00){\makebox(0,0)[cc]{$P_1 $}}
\put(90.00,10.00){\line(1,-1){5.00}}
\put(95.00,5.00){\line(1,1){5.00}}
\put(100.00,10.00){\line(1,-1){5.00}}
\put(105.00,5.00){\line(1,1){5.00}}
\put(110.00,10.00){\line(0,1){10.00}}
\put(110.00,20.00){\line(-1,1){5.00}}
\put(105.00,25.00){\line(-1,-1){5.00}}
\put(100.00,20.00){\line(-1,1){5.00}}
\put(95.00,25.00){\line(-1,-1){5.00}}
\put(90.00,20.00){\line(0,-1){10.00}}
\put(100.00,17.00){\makebox(0,0)[cc]{\footnotesize $Q$}}
\put(100.00,13.00){\makebox(0,0)[cc]{\footnotesize {\em not} symmetric}}
\put(95.00,15.00){\vector(1,0){10.00}}
\put(105.00,5.00){\circle*{2.00}}
\put(95.00,5.00){\circle*{2.00}}
\put(95.00,2.00){\makebox(0,0)[cc]{\footnotesize $a$}}
\put(105.00,2.00){\makebox(0,0)[cc]{\footnotesize $b$}}
\put(107.00,25.00){\makebox(0,0)[lc]{\footnotesize $p_0 = p=\overline{p}_1 $}}
\put(93.00,25.00){\makebox(0,0)[rc]{\footnotesize $d_0 = d=\overline{d}_1 $}}
\put(90.00,30.00){\line(1,-1){5.00}}
\put(95.00,25.00){\line(1,1){5.00}}
\put(100.00,30.00){\line(1,-1){5.00}}
\put(105.00,25.00){\line(1,1){5.00}}
\put(110.00,30.00){\line(0,1){10.00}}
\put(110.00,40.00){\line(-1,1){5.00}}
\put(105.00,45.00){\line(-1,-1){5.00}}
\put(100.00,40.00){\line(-1,1){5.00}}
\put(95.00,45.00){\line(-1,-1){5.00}}
\put(90.00,40.00){\line(0,-1){10.00}}
\put(100.00,37.00){\makebox(0,0)[cc]{\footnotesize $R_1 $}}
\put(100.00,33.00){\makebox(0,0)[cc]{\footnotesize symmetric}}
\put(105.00,25.00){\circle*{2.00}}
\put(95.00,25.00){\circle*{2.00}}
\put(107.00,45.00){\makebox(0,0)[lc]{\footnotesize $p_1 =\overline{p}_2 $}}
\put(93.00,45.00){\makebox(0,0)[rc]{\footnotesize $d_1=\overline{d}_2 $}}
\put(90.00,50.00){\line(1,-1){5.00}}
\put(95.00,45.00){\line(1,1){5.00}}
\put(100.00,50.00){\line(1,-1){5.00}}
\put(105.00,45.00){\line(1,1){5.00}}
\put(110.00,50.00){\line(0,1){10.00}}
\put(110.00,60.00){\line(-1,1){5.00}}
\put(105.00,65.00){\line(-1,-1){5.00}}
\put(100.00,60.00){\line(-1,1){5.00}}
\put(95.00,65.00){\line(-1,-1){5.00}}
\put(90.00,60.00){\line(0,-1){10.00}}
\put(100.00,57.00){\makebox(0,0)[cc]{\footnotesize $R_2 $}}
\put(100.00,53.00){\makebox(0,0)[cc]{\footnotesize symmetric}}
\put(105.00,45.00){\circle*{2.00}}
\put(95.00,45.00){\circle*{2.00}}
\put(107.00,65.00){\makebox(0,0)[lc]{\footnotesize $p_2 =\overline{p}_3 $}}
\put(93.00,65.00){\makebox(0,0)[rc]{\footnotesize $d_2=\overline{d}_3 $}}
\put(105.00,65.00){\circle*{2.00}}
\put(95.00,65.00){\circle*{2.00}}
\put(100.00,75.00){\makebox(0,0)[cc]{$\vdots $}}
\put(107.00,85.00){\makebox(0,0)[lc]{\footnotesize $p_{n-2} =\overline{p}_{n-1} $}}
\put(93.00,85.00){\makebox(0,0)[rc]{\footnotesize $d_{n-2} =\overline{d}_{n-1} $}}
\put(90.00,90.00){\line(1,-1){5.00}}
\put(95.00,85.00){\line(1,1){5.00}}
\put(100.00,90.00){\line(1,-1){5.00}}
\put(105.00,85.00){\line(1,1){5.00}}
\put(110.00,90.00){\line(0,1){10.00}}
\put(110.00,100.00){\line(-1,1){5.00}}
\put(105.00,105.00){\line(-1,-1){5.00}}
\put(100.00,100.00){\line(-1,1){5.00}}
\put(95.00,105.00){\line(-1,-1){5.00}}
\put(90.00,100.00){\line(0,-1){10.00}}
\put(100.00,97.00){\makebox(0,0)[cc]{\footnotesize $R_{n-1} $}}
\put(100.00,93.00){\makebox(0,0)[cc]{\footnotesize symmetric}}
\put(105.00,85.00){\circle*{2.00}}
\put(95.00,85.00){\circle*{2.00}}
\put(107.00,105.00){\makebox(0,0)[lc]{\footnotesize $p_{n-1} =\overline{p}_n $}}
\put(93.00,105.00){\makebox(0,0)[rc]{\footnotesize $d_{n-1} =\overline{d}_n $}}
\put(90.00,110.00){\line(1,-1){5.00}}
\put(95.00,105.00){\line(1,1){5.00}}
\put(100.00,110.00){\line(1,-1){5.00}}
\put(105.00,105.00){\line(1,1){5.00}}
\put(110.00,110.00){\line(0,1){10.00}}
\put(110.00,120.00){\line(-1,1){5.00}}
\put(105.00,125.00){\line(-1,-1){5.00}}
\put(100.00,120.00){\line(-1,1){5.00}}
\put(95.00,125.00){\line(-1,-1){5.00}}
\put(90.00,120.00){\line(0,-1){10.00}}
\put(100.00,117.00){\makebox(0,0)[cc]{\footnotesize $R_n $}}
\put(100.00,113.00){\makebox(0,0)[cc]{\footnotesize symmetric}}
\put(105.00,105.00){\circle*{2.00}}
\put(95.00,105.00){\circle*{2.00}}
\put(107.00,125.00){\makebox(0,0)[lc]{\footnotesize $p_n =\tilde{d} $}}
\put(93.00,125.00){\makebox(0,0)[rc]{\footnotesize $d_n=\tilde{p} $}}
\put(105.00,125.00){\circle*{2.00}}
\put(95.00,125.00){\circle*{2.00}}
\put(90.00,130.00){\line(1,-1){5.00}}
\put(95.00,125.00){\line(1,1){5.00}}
\put(100.00,130.00){\line(1,-1){5.00}}
\put(105.00,125.00){\line(1,1){5.00}}
\put(110.00,130.00){\line(0,1){10.00}}
\put(110.00,140.00){\line(-1,1){5.00}}
\put(105.00,145.00){\line(-1,-1){5.00}}
\put(100.00,140.00){\line(-1,1){5.00}}
\put(95.00,145.00){\line(-1,-1){5.00}}
\put(90.00,140.00){\line(0,-1){10.00}}
\put(100.00,137.67){\makebox(0,0)[cc]{\footnotesize $\tilde{Q} $}}
\put(100.00,133.00){\makebox(0,0)[cc]{\footnotesize {\em not} symmetric}}
\put(107.00,145.00){\makebox(0,0)[lc]{\footnotesize $\tilde{a} $}}
\put(93.00,145.00){\makebox(0,0)[rc]{\footnotesize $\tilde{b} $}}
\put(105.00,145.00){\circle*{2.00}}
\put(95.00,145.00){\circle*{2.00}}
\put(75.00,145.00){\makebox(0,0)[cc]{$P_2 $}}
\put(105.00,135.00){\vector(-1,0){10.00}}
\end{picture}
}

\caption{Ordered sets $P_1 $ and $P_2 $ as constructed in Lemma \protect\ref{severalmiddle}.
The mentioned symmetry is symmetry along the vertical axis.}
\label{switstck}

\end{figure}

\begin{lem}
\label{severalmiddle}

Let $Q$ be an ordered set such that
\begin{enumerate}
\item
$Q$ has exactly two maximal elements $d$ and $p$,

\item
$d$ and $p$ have the same rank,

\item
\label{Qpsi}
There is an isomorphism $\psi :Q\setminus \{ d\} \to Q\setminus \{ p\} $
with $\psi (p)=d$,

\item
$Q$ has two minimal elements $a$ and $b$,

\item
\label{Qpsia}
$Q\setminus \{ a\} $ has an automorphism $\psi ^a $ with $\psi ^a
(p) = d$, $\psi ^a (d) = p$, and $\psi ^a (b) = b$,

\item
\label{Qpsib}
$Q\setminus \{ b\} $ has an automorphism $\psi ^b $ with $\psi ^b
(p) = d$ and $\psi ^b (d) = p$, and $\psi ^b (a) = a$,

\item
$Q$ is rigid.

\setcounter{itemtransfer}{\value{enumi}}

\end{enumerate}

Let $R$ be an ordered set such that
\begin{enumerate}
\setcounter{enumi}{\value{itemtransfer}}

\item
\label{Rtwomax}
$R$ has exactly two maximal elements $d$ and $p$,

\item
\label{Rtwomaxrnk}
$d$ and $p$ have the same rank,

\item
\label{Rtwomin}
$R$ has exactly two minimal elements $\overline{d}$ and
$\overline{p}$,


\item
\label{Rpaut}
$R\setminus \{ \overline{p} \} $ has an automorphism
$\psi ^{\overline{p} } $ such that $\psi ^{\overline{p} }  (d)=p $,
$\psi ^{\overline{p} }  (p)=d $, and $\psi ^{\overline{p} }  (\overline{d} )=\overline{d}  $,

\item
\label{Rdaut}
$R\setminus \{ \overline{d} \} $ has an automorphism
$\psi ^{\overline{d} } $ such that $\psi ^{\overline{d} }  (d)=p $,
$\psi ^{\overline{d} }  (p)=d $,  and $\psi ^{\overline{d} }  (\overline{p} )=\overline{p}  $,

\item
\label{Rswit}
$R$ has an automorphism $\varphi $ with $\varphi (d)=p$, $\varphi (p)=d$,
$\varphi (\overline{d})=\overline{p}$, $\varphi (\overline{p})=\overline{d}$,

\item
\label{notwistR}
$R$ has no automorphism that is the identity on the minimal elements and not the
identity on the maximal elements.

\end{enumerate}

Let $\tilde{Q} $ be the dual of $Q$ and let $R_1 , \ldots , R_n $ be isomorphic copies of
$R$ such that $Q, \tilde{Q} , R_1 , \ldots , R_n $ are all mutually disjoint.
Let the elements of $R_i $ be distinguished by subscripts $i$, that is,
the maximal and minimal elements of $R_i $ are $d_i $, $p_i $ and $\overline{d} _i $, $\overline{p} _i $.
Let the elements of $\tilde{Q}$ similarly be distinguished by tildes.

Define $P_1 $ to be the ordered set obtained from
$Q, R_1 , \ldots , R_n , \tilde{Q} $ as follows (also cf. Figure \ref{switstck}).

\renewcommand{\theenumi}{\roman{enumi}}

\begin{enumerate}
\item
All non-maximal elements of $Q$ are below all non-minimal elements of $R_1 $.

\item
The element $d$ is identified with the element $\overline{d} _1 $ and
the element $p$ is identified with the element $\overline{p} _1 $.
Call the thus obtained elements $d_0 $ and $p_0 $, respectively.

\item
For $i=1, \ldots , n-1$,
all non-maximal elements of $R_i $ are below all non-minimal elements of $R_{i+1} $.

\item
For $i=1, \ldots , n-1$,
the element $d_i $ is identified with the element $\overline{d} _{i+1} $ and
the element $p_i $ is identified with the element $\overline{p} _{i+1} $.
Call the thus obtained elements $d_i $ and $p_i $, respectively.

\item
All non-maximal elements of $R_n $ are below all non-minimal elements of $\tilde{Q} $.

\item
\label{lastmerge}
The element $d_n $ is identified with the element $\tilde{d} $ and
the element $p_n $ is identified with the element $\tilde{p} $.
Call the thus obtained elements $d_n $ and $p_n $, respectively.

\item
Plus all comparabilities forced by transitivity.

\end{enumerate}

Define $P_2 $ to be the ordered set obtained from
$Q, R_1 , \ldots , R_n , \tilde{Q} $ in the same way as $P_1 $ except that
\ref{lastmerge} is replaced with the following (also cf. Figure \ref{switstck}).

\begin{enumerate}
\item[\ref{lastmerge}']
The element $d_n $ is identified with the element $\tilde{p} $ and
the element $p_n $ is identified with the element $\tilde{d} $.
Call the thus obtained elements $d_n $ and $p_n $, respectively.

\end{enumerate}

\renewcommand{\theenumi}{\alph{enumi}}

Then
\begin{enumerate}
\item
\label{P1P2notiso}
$P_1 $ and $P_2 $ are not isomorphic,
\item
\label{P1P2eqmarkedextr}
$P_1 $ and $P_2 $ have
equal marked maximal and minimal decks,

\item
\label{dipieqmark}
For $i=0, \ldots , n$
the
card $P_1 \setminus \{ d_i \} $ is isomorphic to the card $P_2\setminus \{ d_i \} $ and
the
card $P_1 \setminus \{ p_i \} $ is isomorphic to the card $P_2\setminus \{ p_i \} $,
and for all isomorphisms $\Phi $ of cards and all
elements $x\in P_1 $ we have that ${\rm rank} _{P_2 } (\Phi (x) )={\rm rank} _{P_1 } (x)$.

\item
\label{P1P2eqnremhood}
For all the above mentioned isomorphic cards $P_1 \setminus \{ x\} $
and $P_2 \setminus \{ x\} $ the
neighborhoods $\updownarrow _{P_1 } x $ and $\updownarrow _{P_2 } x $
are isomorphic.

\end{enumerate}

\renewcommand{\theenumi}{\arabic{enumi}}

\end{lem}

{\bf Proof.}
To prove \ref{P1P2notiso}, we assume that $P_1 $ and $P_2 $ are isomorphic.
So suppose that $\Phi :P_1 \to P_2 $ is an isomorphism.
Then, because $Q$ is rigid, we have that $\Phi (a)=a$, $\Phi (b)=b$,
$\Phi (d_0 ) = d_0 $ and
$\Phi (p_0 )= p_0 $.
By property \ref{notwistR},
this implies that
$\Phi (d_1 ) = d_1 $,
$\Phi (p_1 )= p_1 $, $\ldots $, $\Phi (d_n ) = d_n $,
$\Phi (p_n )= p_n $.
But then, because in $P_1 $ we have $d_n = \tilde{d} $ and $p_n = \tilde{p}$, while
in $P_2 $ we have $d_n = \tilde{p} $ and $p_n = \tilde{d}$, $\Phi |_{\tilde{Q} } $ would be an
automorphism of $\tilde{Q}$ with
$
{ \Phi |_{\tilde{Q} } ( \tilde{d} ) = \tilde{p} } $ and
${ \Phi |_{\tilde{Q} } ( \tilde{p} ) = \tilde{d} } $.
This is a contradiction to the rigidity of $\tilde{Q} $.
Therefore, $P_1 $ and $P_2 $ cannot be isomorphic.

To show
\ref{P1P2eqmarkedextr},
first note that
$P_1 \setminus \{ a\} $ is isomorphic to $P_2 \setminus \{ a\} $.
To see this
let $\varphi _i $ denote the automorphism for $R_i $ guaranteed by
property \ref{Rswit}.
We define

$$\Phi (x):=\cases{ \psi ^a (x); & if $x\in Q\setminus \{ a\} $, \cr
\varphi _i (x); & if $x\in R_i $, \cr
x; & if $x\in \tilde{Q}\setminus \{ \tilde{d}, \tilde{p} \} $. \cr
} $$

The function $\Phi $ is well-defined and bijective between
$P_1 \setminus \{ a\} $ and $P_2 \setminus \{ a\} $ and it maps minimal elements of $P_1 $
to minimal elements of $P_2 $.
To see that $\Phi $ is order-preserving both ways, let $x<y$ in $P_1$.
It is trivial to infer that $\Phi (x)<\Phi (y)$ is equivalent to $x<y$ unless
$x\in \{ \tilde{d} , \tilde{p} \} $.
Assume without loss of generality that $x=\tilde{d} = d_n $.
Then $\Phi (x) =p_n = \tilde{d} =x$.
Since $y>x$ we have $y\in \tilde{Q} $ and thus $\Phi (y)=y$. The other direction, as well as the proof for
$x=\tilde{p} = p_n $ is similar.

We have shown that
$P_1 \setminus \{ a\} $ is isomorphic to $P_2 \setminus \{ a\} $ and that the isomorphism
preserves the marked property ``minimality".
Similarly, $P_1 \setminus \{ b\} $ is isomorphic to $P_2 \setminus \{ b\} $
(and minimal elements of $P_1 $ are mapped to minimal elements of $P_2 $) via

$$\Phi (x):=\cases{ \psi ^b (x); & if $x\in Q\setminus \{ b\} $, \cr
\varphi _i (x); & if $x\in R_i $, \cr
x; & if $x\in \tilde{Q}\setminus \{ \tilde{d}, \tilde{p} \} $. \cr
} $$

The proof that $P_1 $ and $P_2 $ have equal marked maximal decks is similar.
The set
$P_1 \setminus \{ \tilde{a} \} $ is isomorphic to $P_2 \setminus \{ \tilde{a} \} $ via

$$\Phi (x):=\cases{ \psi ^{\tilde{a} } (x); & if $x\in \tilde{Q}\setminus \{ \tilde{a}\} $, \cr
x; & if $x\not\in \tilde{Q} $,  \cr
} $$

where $\psi ^{\tilde{a} } $ denotes the automorphism of $\tilde{Q} \setminus \{ \tilde{a} \} $
that is guaranteed by the dual of property \ref{Qpsia}.
Clearly $\Phi $ maps maximal elements of $P_1 $ to maximal elements of $P_2 $.
The set
$P_1 \setminus \{ \tilde{b} \} $ is isomorphic to $P_2 \setminus \{ \tilde{b} \} $
(and maximal elements of $P_1 $ are mapped to maximal elements of $P_2 $)
via

$$\Phi (x):=\cases{ \psi ^{\tilde{b} } (x); & if $x\in \tilde{Q}\setminus \{ \tilde{b}\} $, \cr
x; & if $x\not\in \tilde{Q} $,  \cr
} $$

where $\psi ^{\tilde{b} } $ denotes the automorphism of $\tilde{Q} \setminus \{ \tilde{b} \} $
that is guaranteed by the dual of property \ref{Qpsib}.

In regards to \ref{P1P2eqnremhood} note that the above isomorphisms show that
for $x\in \{ a,b,\tilde{a},\tilde{b} \} $
the
neighborhoods $\updownarrow _{P_1 } x $ and $\updownarrow _{P_2 } x $
are isomorphic.
(For example, the isomorphism between $P_1 \setminus \{ b\} $ and
$P_2 \setminus \{ b\} $ provides an isomorphism
between
$\updownarrow _{P_1 } a= \uparrow _{P_1 } a$ and
$\updownarrow _{P_2 } a =\uparrow _{P_2 } a$.)

For \ref{dipieqmark}, let it be stated here that it is easy to see that
all isomorphisms $\Phi $ constructed in the following satisfy
${\rm rank} _{P_2 } (\Phi (x) )={\rm rank} _{P_1 } (x)$ for all
elements $x\in P_1 $.

Now first notice that
the set
$P_1 \setminus \{ d_n \} $ is isomorphic to $P_2 \setminus \{ d_n \} $ via
$$\Phi (x):=\cases{ \tilde{\psi }  (x); & if $x\in \tilde{Q}\setminus \{ \tilde{d}\} $, \cr
x; & if $x\not\in \tilde{Q} $,   \cr
} $$

where $\tilde{\psi } $ is the isomorphism guaranteed by the dual of
property \ref{Qpsi}.
The set
$P_1 \setminus \{ p_n \} $ is isomorphic to $P_2 \setminus \{ p_n \} $ via
$$\Phi (x):=\cases{ \left( \tilde{\psi } \right) ^{-1} (x); & if $x\in \tilde{Q}\setminus \{ \tilde{p}\} $, \cr
x; & if $x\not\in \tilde{Q} $.   \cr
} $$

Finally let $i\in \{ 0,\ldots , n-1\} $.
For $R_i $ denote the automorphisms $\psi ^{\overline{p} }  $ and $\psi ^{\overline{d} }  $
of the respective cards of $R$
guaranteed by properties
\ref{Rpaut} and \ref{Rdaut} by $\psi ^{\overline{p} }  _i $ and $\psi ^{\overline{d} }  _i $.
Then
the set
$P_1 \setminus \{ p_i \} $ is isomorphic to $P_2 \setminus \{ p_i \} $ via

$$\Phi (x):=\cases{
x; & if $x\in \tilde{Q}\setminus \{ \tilde{d}, \tilde{p} \} $ or $x\in R_j $ for $j\leq i$ or $x\in Q$, \cr
\varphi _j (x); & if $x\in R_j $ for $j>i+1$, \cr
\psi ^{\overline{p} }  _i (x) ; & if $x\in R_{i+1} \setminus \{ p_i \} $.\cr
} $$

All parts of the definition of isomorphism are readily verified.
Similarly,
the set
$P_1 \setminus \{ d_i \} $ is isomorphic to $P_2 \setminus \{ d_i \} $ via

$$\Phi (x):=\cases{
x; & if $x\in \tilde{Q}\setminus \{ \tilde{d}, \tilde{p} \} $ or $x\in R_j $ for $j\leq i$ or $x\in Q$, \cr
\varphi _j (x); & if $x\in R_j $ for $j>i+1$, \cr
\psi ^{\overline{d} }  _i (x) ; & if $x\in R_{i+1} \setminus \{ d_i \} $.\cr
} $$

To finish the proof of \ref{P1P2eqnremhood}, similar to
what was said after the proof of \ref{P1P2eqmarkedextr},
the above isomorphisms show that
for $x\in \{ d_0 , p_0 , \ldots , d_{n-1} , p_{n-1} \} $
the
neighborhoods $\updownarrow _{P_1 } x $ and $\updownarrow _{P_2 } x $
are isomorphic.
Just use the isomorphism between the cards on which the respective
``other" element of that rank has been erased.
\qed

\begin{figure}

\centerline{
\unitlength 1mm 
\linethickness{0.4pt}
\ifx\plotpoint\undefined\newsavebox{\plotpoint}\fi 
\begin{picture}(146,93)(0,0)
\put(50,50){\line(1,-4){2.5}}
\put(15,40){\line(1,4){2.5}}
\put(60,50){\line(1,-4){2.5}}
\put(25,40){\line(1,4){2.5}}
\put(70,50){\line(1,-4){2.5}}
\put(35,40){\line(1,4){2.5}}
\put(85,50){\line(1,-4){2.5}}
\put(120,40){\line(1,4){2.5}}
\put(95,50){\line(1,-4){2.5}}
\put(130,40){\line(1,4){2.5}}
\put(105,50){\line(1,-4){2.5}}
\put(140,40){\line(1,4){2.5}}
\put(52.5,40){\line(1,4){2.5}}
\put(17.5,50){\line(1,-4){2.5}}
\put(62.5,40){\line(1,4){2.5}}
\put(27.5,50){\line(1,-4){2.5}}
\put(72.5,40){\line(1,4){2.5}}
\put(37.5,50){\line(1,-4){2.5}}
\put(87.5,40){\line(1,4){2.5}}
\put(122.5,50){\line(1,-4){2.5}}
\put(97.5,40){\line(1,4){2.5}}
\put(132.5,50){\line(1,-4){2.5}}
\put(107.5,40){\line(1,4){2.5}}
\put(142.5,50){\line(1,-4){2.5}}
\put(50,50){\circle*{1}}
\put(15,40){\circle*{1}}
\put(60,50){\circle*{1}}
\put(25,40){\circle*{1}}
\put(70,50){\circle*{1}}
\put(35,40){\circle*{1}}
\put(85,50){\circle*{1}}
\put(120,40){\circle*{1}}
\put(95,50){\circle*{1}}
\put(130,40){\circle*{1}}
\put(105,50){\circle*{1}}
\put(140,40){\circle*{1}}
\put(52.5,40){\circle*{1}}
\put(17.5,50){\circle*{1}}
\put(62.5,40){\circle*{1}}
\put(27.5,50){\circle*{1}}
\put(72.5,40){\circle*{1}}
\put(37.5,50){\circle*{1}}
\put(87.5,40){\circle*{1}}
\put(122.5,50){\circle*{1}}
\put(97.5,40){\circle*{1}}
\put(132.5,50){\circle*{1}}
\put(107.5,40){\circle*{1}}
\put(142.5,50){\circle*{1}}
\put(55,50){\circle*{1}}
\put(20,40){\circle*{1}}
\put(65,50){\circle*{1}}
\put(30,40){\circle*{1}}
\put(75,50){\circle*{1}}
\put(40,40){\circle*{1}}
\put(90,50){\circle*{1}}
\put(125,40){\circle*{1}}
\put(100,50){\circle*{1}}
\put(135,40){\circle*{1}}
\put(110,50){\circle*{1}}
\put(145,40){\circle*{1}}
\put(75,50){\line(1,-4){5}}
\put(80,30){\line(1,4){5}}
\put(80,30){\circle*{2}}
\qbezier(55,50)(57.5,31.5)(80,30)
\qbezier(105,50)(102.5,31.5)(80,30)
\qbezier(80,30)(67.5,35)(65,50)
\qbezier(80,30)(92.5,35)(95,50)
\put(75,50){\circle{2}}
\put(40,40){\circle{2}}
\put(60,50){\circle{2}}
\put(25,40){\circle{2}}
\put(90,50){\circle{2}}
\put(125,40){\circle{2}}
\put(95,50){\circle{2}}
\put(130,40){\circle{2}}
\put(69,49){\framebox(2,2)[cc]{}}
\put(34,39){\framebox(2,2)[]{}}
\put(54,49){\framebox(2,2)[cc]{}}
\put(19,39){\framebox(2,2)[]{}}
\put(84,49){\framebox(2,2)[cc]{}}
\put(119,39){\framebox(2,2)[]{}}
\put(109,49){\framebox(2,2)[cc]{}}
\put(144,39){\framebox(2,2)[]{}}
\put(80,70){\circle*{2}}
\qbezier(20,40)(23,66)(80,70)
\qbezier(140,40)(137,66)(80,70)
\qbezier(80,70)(33,64)(30,40)
\qbezier(80,70)(127,64)(130,40)
\qbezier(80,70)(44.5,62)(40,40)
\qbezier(80,70)(115.5,62)(120,40)
\bezier{40}(15,49)(27.5,47)(40,49)
\bezier{40}(145,49)(132.5,47)(120,49)
\bezier{40}(75,41)(62.5,43)(50,41)
\bezier{40}(85,41)(97.5,43)(110,41)
\put(55,41.75){\line(-3,1){20}}
\put(105,41.75){\line(3,1){20}}
\put(83,27){\makebox(0,0)[cc]{$c_b $}}
\put(83,73){\makebox(0,0)[cc]{$c_t $}}
\put(17.5,37){\makebox(0,0)[cc]{$A_1 $}}
\put(52.5,53){\makebox(0,0)[cc]{$B_1 $}}
\put(97.5,53){\makebox(0,0)[cc]{$C_1 $}}
\put(132.5,37){\makebox(0,0)[cc]{$D_1 $}}
\put(27.5,37){\makebox(0,0)[cc]{$A_2 $}}
\put(62.5,53){\makebox(0,0)[cc]{$B_2 $}}
\put(107.5,53){\makebox(0,0)[cc]{$C_2 $}}
\put(142.5,37){\makebox(0,0)[cc]{$D_2 $}}
\put(37.5,37){\makebox(0,0)[cc]{$A_{1,2} $}}
\put(72.5,53){\makebox(0,0)[cc]{$B_{1,2} $}}
\put(87.5,53){\makebox(0,0)[cc]{$C_{1,2} $}}
\put(122.5,37){\makebox(0,0)[cc]{$D_{1,2} $}}
\put(75,10){\circle*{2}}
\put(75,90){\circle*{2}}
\put(85,10){\circle*{2}}
\put(85,90){\circle*{2}}
\put(75,7){\makebox(0,0)[cc]{$\overline{d}$}}
\put(75,93){\makebox(0,0)[cc]{${d}$}}
\put(85,7){\makebox(0,0)[cc]{$\overline{p}$}}
\put(85,93){\makebox(0,0)[cc]{${p}$}}
\put(15,85){\makebox(0,0)[lc]{$A=A_1 \cup A_2 \cup A_{1,2} $}}
\put(15,75){\makebox(0,0)[lc]{$B=B_1 \cup B_2 \cup B_{1,2} $}}
\put(145,85){\makebox(0,0)[rc]{$C=C_1 \cup C_2 \cup C_{1,2} $}}
\put(145,75){\makebox(0,0)[rc]{$D=D_1 \cup D_2 \cup D_{1,2} $}}
\end{picture}
}

\caption{An ordered set $R$ as needed in Lemma \protect\ref{severalmiddle} and constructed in
Lemma \ref{getsmallerR}.
The middle levels form an ordered set of height $1$.
The maximal element $d$ is above all maximal elements of the
middle levels {\em except} the circled maximal elements
of the middle levels.
The maximal element $p$ is above all maximal elements of the
middle levels {\em except} the boxed maximal elements
of the middle levels.
Similarly,
the minimal element $\overline{d}$ is below all minimal elements of the
middle levels {\em except} the circled minimal elements and
the minimal element $\overline{p}$ is below all minimal elements of the
middle levels {\em except} the boxed minimal elements.
Within the middle levels, the connected dotted arches indicate that the
elements immediately above the upper arch and the elements immediately below the lower arch form a complete
bipartite.
}
\label{flatswit}

\end{figure}

\begin{lem}
\label{getsmallerR}

There is
an ordered set $R$ as described in Lemma \ref{severalmiddle}.

\end{lem}

{\bf Proof.}
Let $R$ be the ordered set indicated in Figure \ref{flatswit}.
We claim that $R$ is a set as desired in the description of the set $R$ in
Lemma \ref{severalmiddle}.

Claims \ref{Rtwomax}, \ref{Rtwomaxrnk}, \ref{Rtwomin} are trivial.
The rest of the proof relies on the following, easy to verify, properties
of the sets $A$, $B$, $C$ and $D$.
First, $A\cup \{ c_t , \overline{d} ,\overline{p} \} $ is dually isomorphic to
$B\cup \{ c_b , d,p\} $ and
$D\cup \{ c_t , \overline{d} ,\overline{p} \} $ is dually isomorphic to
$C\cup \{ c_b , d,p\} $.
Therefore, properties proved for $B$ and $C$ will hold dually for $A$ and $D$.
The set $B\cup \{ c_b , d,p\} $ is rigid.
This is because any automorphism of $B\cup \{ c_b , d,p\} $ must map $c_b $ to itself
and the sets $B_1 $, $B_2 $ and $B_{1,2} $ to themselves, respectively.
This implies that the automorphism must be the identity on $B$ and then it must be the identity on
$\{ d,p\} $ also.

Moreover,
there is exactly one isomorphism
from
$B\cup \{ c_b , d,p\} $ to $C\cup \{ c_b , d,p\} $. This isomorphism
maps $d$ to $p$, $p$ to $d$, $c_b $ to $c_b $ and $B_{1,2} $ to $C_{1,2} $, $B_{1} $ to $C_{1} $ and
$B_{2} $ to $C_{2} $.
Furthermore, there is exactly one isomorphism from
$B\cup \{ c_b , d\} $ to $C\cup \{ c_b , d\} $.
This isomorphism maps $d$ to $d$, $c_b $ to $c_b $,
$B_{1,2} $ to $C_{1} $, $B_{2} $ to $C_{1,2} $ and $B_{1} $ to $C_{2} $.
Similarly, there is exactly one isomorphism
from
$B\cup \{ c_b , p\} $ to $C\cup \{ c_b , p\} $.
This isomorphism maps $p$ to $p$, $c_b $ to $c_b $,
$B_{1,2} $ to $C_{2} $, $B_{2} $ to $C_{1} $ and $B_{1} $ to $C_{1,2} $.
These facts and their duals will be used freely in the following.

For Claim \ref{Rpaut} define $\psi ^{\overline{p} } $ to be
\begin{enumerate}
\item
$\psi ^{\overline{p} } (\overline{d} ):=\overline{d} $,
$\psi ^{\overline{p} } (d ):=p $,
$\psi ^{\overline{p} } (p ):=d $,

\item
On $A\cup \{ c_t \} $, $\psi ^{\overline{p} } $ is
the restriction of
the unique isomorphism from
$A\cup \{ c_t , \overline{d} \} $ to
$D\cup \{ c_t , \overline{d} \} $,

\item
On $D\cup \{ c_t \} $, $\psi ^{\overline{p} } $ is
the restriction of
the unique isomorphism from
$D\cup \{ c_t , \overline{d} \} $ to
$A\cup \{ c_t \overline{d} \} $,

\item
On $B\cup \{ c_b \} $, $\psi ^{\overline{p} } $ is
the restriction of the unique
isomorphism from
$B\cup \{ c_b, d,p \} $ to
$C\cup \{ c_b, d,p \} $,

\item
On $C\cup \{ c_b \} $, $\psi ^{\overline{p} } $ is
the restriction of the unique
isomorphism from
$C\cup \{ c_b, d,p \} $ to
$B\cup \{ c_b, d,p \} $.

\end{enumerate}

Then $\psi ^{\overline{p} } $ is as desired.
Claim \ref{Rdaut} is proved similarly.

For Claim \ref{Rswit} define $\varphi $ to be
\begin{enumerate}
\item
$\varphi (\overline{d} ):=\overline{p} $,
$\varphi (\overline{p} ):=\overline{d} $,
$\varphi (d ):=p $,
$\varphi (p ):=d $,

\item
On $A\cup \{ c_t \} $, $\varphi $ is
the restriction of the unique
isomorphism from
$A\cup \{ c_t , \overline{d}, \overline{p} \} $ to
$D\cup \{ c_t, \overline{d}, \overline{p} \} $,

\item
On $D\cup \{ c_t \} $, $\varphi $ is
the restriction of the unique
isomorphism from
$D\cup \{ c_t, \overline{d}, \overline{p} \} $ to
$A\cup \{ c_t \overline{d}, \overline{p} \} $,

\item
On $B\cup \{ c_b\} $, $\varphi $ is
the restriction of the unique
isomorphism from
$B\cup \{ c_b, d,p \} $ to
$C\cup \{ c_b , d,p \} $,

\item
On $C\cup \{ c_b\} $, $\varphi $ is
the restriction of the unique
isomorphism from
$C\cup \{ c_b , d,p \} $ to
$B\cup \{ c_b, d,p \} $.

\end{enumerate}

Finally, for Claim \ref{notwistR} suppose the automorphism
$\Psi :R\to R$
is the identity on the minimal elements.
Then, because the unique isomorphism between
$A\cup \{ c_t, \overline{d},\overline{p} \} $ and $D\cup \{ c_t ,\overline{d},\overline{p} \} $
switches $\overline{d} $ and $\overline{p} $,
$\Psi $ must map $A \cup \{c_t \} $ to $A \cup \{ c_t\} $ and
$D  \cup \{ c_t\} $ to $D \cup \{ c_t\}  $.
Therefore, $\Psi $ must map
$B \cup \{ c_b\} $ to $B \cup \{ c_b\}  $ and
$C \cup \{ c_b\}  $ to $C \cup \{ c_b\}  $.
Since $B\cup \{ c_b , d , p \} $ is
rigid, this means that
$\Psi $ must fix $d $ and $p $.
\qed

\vspace{.1in}

The set $R$ in Figure \ref{flatswit} has an additional property that will allow us to
prove further properties of our examples.

\begin{lem}
\label{takeoutc}

The sets $R\setminus \{ c_b\} $ and $R\setminus \{ c_t \} $
with $R$ as in Figure \ref{flatswit} each have an automorphism
$\Psi $ with $\Psi (\overline{d} )=\overline{d} $, $\Psi (\overline{p} )=\overline{p} $,
$\Psi (d)=p$ and $\Psi (p)=d$.

\end{lem}

{\bf Proof.}
For $R\setminus \{ c_t \} $ we define $\Psi (\overline{d} ):=\overline{d} $, $\Psi (\overline{p} ):=\overline{p} $,
$\Psi (d):=p$ and $\Psi (p):=d$. We let
$\Psi $ map $A_{1,2} $ to $D_{1,2} $, $A_{2} $ to $D_{1} $, $A_{1} $ to $D_{2} $ and vice versa.
(Visually, in each case the map is obtained by sliding one wedge horizontally onto the other.)

On $B\cup \{ c_b, d,p\} $ we define $\Psi $ to be the unique isomorphism from
$B\cup \{ c_b, d,p\} $ to $C\cup \{ c_b, d,p\} $.
Finally, on $C\cup \{ c_b, d,p\} $ we define $\Psi $ to be the unique isomorphism from
$C\cup \{ c_b, d,p\} $ to $B\cup \{ c_b, d,p\} $.

For $R\setminus \{ c_b \} $ we let $\Psi $ be the identity on $\{ \overline{d} , \overline{p} , c_t \} \cup A\cup D$.
We let $\Psi (d)=p$ and $\Psi (p)=d$. Finally $\Psi $ maps $B_{1,2} $ to itself, $B_1 $ to $B_2 $, $B_2 $ to $B_1 $,
$C_{1,2} $ to itself, $C_1 $ to $C_2 $, $C_2 $ to $C_1 $, each in such a way that the comparabilities
with the appropriate maximal elements are preserved.
\qed

\vspace{.1in}

With Lemmas \ref{severalmiddle}
and \ref{getsmallerR} proved, we can state our first main result.
Aside from insights on ranks and maximal and minimal cards, we also see that
our construction yields pairs of nonisomorphic sets
for which a significant
number of cards is isomorphic.

\begin{define}

For two ordered sets $P_1 $ and $P_2 $ with $n$ elements each, define the
{\bf equal card ratio} $ECR(P_1 , P_2 )$ to be the number of
isomorphic cards divided by the size of the set.

\end{define}

\begin{theorem}
\label{tallones}

There is a sequence of pairs of ordered sets $(P_1 ^n , P_2 ^n )$
such that
\begin{enumerate}
\item
$P_1 ^n $ is not isomorphic to $P_2 ^n $,
\item
\label{eqmaxmindeck}
$P_1 ^n $ and $P_2 ^n $ have equal marked maximal and minimal decks,
\item
\label{eqrankkdecks}
There are ranks $k_0, \ldots , k_n $ such that $P_1 ^n $ and $P_2 ^n $ have
equal marked rank $k_i $ decks,
\item
For any
pair of isomorphic cards as in parts \ref{eqmaxmindeck} and \ref{eqrankkdecks},
the neighborhoods of the respective removed elements are isomorphic,
\item
\label{onlyfourleft}
There are only four ranks that do not produce any isomorphic cards,
\item
\label{ecr>0}
$\displaystyle{ \liminf _{n\to \infty } ECR(P_1 ^n , P_2 ^n )\geq {1\over 10} } $,
\item
\label{equalnhooddecks}
For every rank $k$, the rank $k$ neighborhood decks of $P_1 ^n $ and $P_2 ^n $
are equal.
\end{enumerate}

\end{theorem}

{\bf Proof.}
To construct sets as indicated,
with notation as in Lemma \ref{severalmiddle}, we make the following choices.
\begin{enumerate}
\item
For all $n$, as the set $Q$, use a fixed set $Q$ as in the proof of Theorem 5.3 of \cite{Schmoreexam}.
These sets have height $3$.
\item
For all $n$, as the set $R$, use a fixed set $R$ as guaranteed by
Lemma \ref{getsmallerR}.
\item
The ranks $k_0 , \ldots , k_n $ are
$k_i = {\rm height} (Q)+i\cdot {\rm height} (R)$.
\end{enumerate}

This construction, independent of the choices for $Q$ and $R$,
yields a sequence of sets that satisfy all parts of this theorem except
possibly parts \ref{onlyfourleft}, \ref{ecr>0} and \ref{equalnhooddecks}.
For these parts the construction must be done with the set $R$ indicated in the proof of
Lemma \ref{getsmallerR}.

For part \ref{onlyfourleft} note that
by Lemma \ref{takeoutc} the erasure of an element $c_t $ or $c_b $
at corresponding ranks produces isomorphic marked cards for $P_1 ^n $ and $P_2 ^n $.
This means that only at the ranks $1$ and $2$ and at the two ranks immediately below the
maximal elements will there be no elements whose removal produces isomorphic cards.

For part \ref{ecr>0},
first note that because each pair of sets has at least $4n+6$ equal cards we have
$$\displaystyle{ ECR(P_1 ^n , P_2 ^n )\geq {4n+6\over n(|R|-2)+2|Q|-2}
\ontop{n\to \infty }{\longrightarrow } {4\over |R|-2}. } $$

That is,
our lower bound on
$\displaystyle{
\liminf _{n\to \infty } ECR(P_1 ^n , P_2 ^n )} $ will solely depend on the number of elements that $R$ has.

With the set $R$ as given in the proof of Lemma \ref{getsmallerR} we have $|R|=42$ and so
$\displaystyle{ \liminf _{n\to \infty } ECR(P_1 ^n , P_2 ^n )\geq {4\over 42-2}={1\over 10} .} $

For part \ref{equalnhooddecks} we know that
the neighborhoods $\updownarrow _{P_1 ^n } x$ and $\updownarrow _{P_2 ^n } x$
are isomorphic for $x\in \{ a,b, \tilde{a}, \tilde{b},d_0 , p_0 , \ldots , d_n , p_n \} $.
This leaves the non-extremal elements of the $R_i $ and the non-extremal
elements of $Q$ and $\tilde{Q} $.

We first consider the non-extremal elements of the $R_i $.
For $x\in A\cup D$
we have that $\uparrow _R x\setminus \{ x\} $ is a four-crown if $x$ is minimal and
a two-antichain if $x$ is maximal.
Similarly,
$\uparrow _R x\setminus \{ x\} $ is
a two-antichain for the maximal elements of $C$ and $D$ that are below both
$d$ and $p$ as well as for $c_t $.
This means that $\updownarrow _R x$ has an automorphism that fixes
$\overline{d} $ and $\overline{p} $ and that switches $d$ and $p$,
which means that for these elements
(in any of the $R_i $) we have that
$\updownarrow _{P_1 ^n } x$ and $\updownarrow _{P_2 ^n } x$
are isomorphic.
For the minimal elements of $B_{1,2} $ and $C_{1,2} $, call them
$b_{1,2} $ and $c_{1,2} $, the set
$\uparrow _R x\setminus \{ x\} $ is the disjoint union of two 2-chains
and for these elements (in any $R_i $) we have
$\updownarrow _{P_1 ^n } b_{1,2} $ is isomorphic to
$\updownarrow _{P_2 ^n } c_{1,2} $
and
$\updownarrow _{P_1 ^n } c_{1,2} $ is isomorphic to
$\updownarrow _{P_2 ^n } b_{1,2} $.
For the minimal elements of $B_1 $ and $B_2 $, call them
$b_1 $ and $b_2 $,
the set
$\uparrow _R x\setminus \{ x\} $ is an ``N"
and for these elements (in any $R_i $) we have
$\updownarrow _{P_1 ^n } b_{1} $ is isomorphic to
$\updownarrow _{P_2 ^n } b_{2} $.
The corresponding elements of $C$ are handled similarly.
The remaining maximal elements of $B$ and $C$ are below exactly one of $d$ and $p$,
call this element $m$.
For these elements (in any $R_i $) we prove that
$\updownarrow _{P_1 ^n } m $ is isomorphic to
$\updownarrow _{P_2 ^n } m $ as follows.
In an $R_i $ with $i<n$ this is because
$R_{i+1} \setminus \{ m\} $ has an automorphism that
switches the maximal elements of $R_{i+1} $.
In $R_n $ this is because $\tilde{Q} \setminus \{ \tilde{d} \} $ is
isomorphic to $\tilde{Q} \setminus \{ \tilde{p} \} $.
Finally,
$\uparrow _R c_b \setminus \{ c_b \} $ is
made up of two sets ``M" with their maximal elements identified.
This means that
$\updownarrow _{P_1 ^n } c_b $ is isomorphic to
$\updownarrow _{P_2 ^n } c_b $.

To consider the non-extremal elements of $Q$, we
need to use sets $Q$ as in
\cite{Schmoreexam}.
For all non-extremal elements $x$ of $Q$ the same ideas as above show that
$\updownarrow _{P_1 ^n } x$ and $\updownarrow _{P_2 ^n } x$
are isomorphic.
If the strict upper bounds of $x$ in $Q$ are a set
(such as a two antichain, a four crown or the disjoint union of
two $2$-chains)
that allows
interchanging the maximal elements, we get an isomorphism.
If not, the strict upper bounds in $Q$ are a singleton, in which
case the switch can be done in $R_1 $, or the strict upper bounds
form an ``N". By choosing the appropriate
permutations in the construction of the set $Q$, this last situation
can be avoided. This takes care of the non-extremal elements of $Q$.
The argument for the non-extremal elements of $\tilde{Q} $ is the dual
of the above.
\qed

%
%

\begin{figure}

\centerline{
\unitlength 1.00mm 
\linethickness{0.4pt}
\ifx\plotpoint\undefined\newsavebox{\plotpoint}\fi 
\begin{picture}(146,93)(0,0)
\put(50,50){\line(1,-4){2.5}}
\put(15,50){\line(1,-4){2.5}}
\put(60,50){\line(1,-4){2.5}}
\put(25,50){\line(1,-4){2.5}}
\put(70,50){\line(1,-4){2.5}}
\put(35,50){\line(1,-4){2.5}}
\put(85,50){\line(1,-4){2.5}}
\put(120,50){\line(1,-4){2.5}}
\put(95,50){\line(1,-4){2.5}}
\put(130,50){\line(1,-4){2.5}}
\put(105,50){\line(1,-4){2.5}}
\put(140,50){\line(1,-4){2.5}}
\put(52.5,40){\line(1,4){2.5}}
\put(17.5,40){\line(1,4){2.5}}
\put(62.5,40){\line(1,4){2.5}}
\put(27.5,40){\line(1,4){2.5}}
\put(72.5,40){\line(1,4){2.5}}
\put(37.5,40){\line(1,4){2.5}}
\put(87.5,40){\line(1,4){2.5}}
\put(122.5,40){\line(1,4){2.5}}
\put(97.5,40){\line(1,4){2.5}}
\put(132.5,40){\line(1,4){2.5}}
\put(107.5,40){\line(1,4){2.5}}
\put(142.5,40){\line(1,4){2.5}}
\put(50,50){\circle*{1}}
\put(15,50){\circle*{1}}
\put(60,50){\circle*{1}}
\put(25,50){\circle*{1}}
\put(70,50){\circle*{1}}
\put(35,50){\circle*{1}}
\put(85,50){\circle*{1}}
\put(120,50){\circle*{1}}
\put(95,50){\circle*{1}}
\put(130,50){\circle*{1}}
\put(105,50){\circle*{1}}
\put(140,50){\circle*{1}}
\put(52.5,40){\circle*{1}}
\put(17.5,40){\circle*{1}}
\put(62.5,40){\circle*{1}}
\put(27.5,40){\circle*{1}}
\put(72.5,40){\circle*{1}}
\put(37.5,40){\circle*{1}}
\put(87.5,40){\circle*{1}}
\put(122.5,40){\circle*{1}}
\put(97.5,40){\circle*{1}}
\put(132.5,40){\circle*{1}}
\put(107.5,40){\circle*{1}}
\put(142.5,40){\circle*{1}}
\put(55,50){\circle*{1}}
\put(20,50){\circle*{1}}
\put(65,50){\circle*{1}}
\put(30,50){\circle*{1}}
\put(75,50){\circle*{1}}
\put(40,50){\circle*{1}}
\put(90,50){\circle*{1}}
\put(125,50){\circle*{1}}
\put(100,50){\circle*{1}}
\put(135,50){\circle*{1}}
\put(110,50){\circle*{1}}
\put(145,50){\circle*{1}}
\put(75,50){\line(1,-4){5}}
\put(80,30){\line(1,4){5}}
\put(80,30){\circle*{2}}
\qbezier(55,50)(57.5,31.5)(80,30)
\qbezier(105,50)(102.5,31.5)(80,30)
\qbezier(80,30)(67.5,35)(65,50)
\qbezier(80,30)(92.5,35)(95,50)
\put(75,50){\circle{2}}
\put(40,50){\circle{2}}
\put(60,50){\circle{2}}
\put(25,50){\circle{2}}
\put(90,50){\circle{2}}
\put(125,50){\circle{2}}
\put(95,50){\circle{2}}
\put(130,50){\circle{2}}
\put(69,49){\framebox(2,2)[cc]{}}
\put(34,49){\framebox(2,2)[]{}}
\put(54,49){\framebox(2,2)[cc]{}}
\put(19,49){\framebox(2,2)[]{}}
\put(84,49){\framebox(2,2)[cc]{}}
\put(119,49){\framebox(2,2)[]{}}
\put(109,49){\framebox(2,2)[cc]{}}
\put(144,49){\framebox(2,2)[]{}}
\put(80,20){\circle*{2}}
\qbezier(20,50)(23,24)(80,20)
\qbezier(140,50)(137,24)(80,20)
\qbezier(80,20)(33,26)(30,50)
\qbezier(80,20)(127,26)(130,50)
\qbezier(80,20)(44.5,28)(40,50)
\qbezier(80,20)(115.5,28)(120,50)
\bezier{40}(15,49)(27.5,47)(40,49)
\bezier{40}(145,49)(132.5,47)(120,49)
\bezier{40}(75,41)(62.5,43)(50,41)
\bezier{40}(85,41)(97.5,43)(110,41)
\put(55,41.75){\line(-3,1){20}}
\put(105,41.75){\line(3,1){20}}
\put(83,27){\makebox(0,0)[cc]{$c_b $}}
\put(83,17){\makebox(0,0)[]{$c_t $}}
\put(17.5,53){\makebox(0,0)[]{$A_1 $}}
\put(52.5,53){\makebox(0,0)[cc]{$B_1 $}}
\put(97.5,53){\makebox(0,0)[cc]{$C_1 $}}
\put(132.5,53){\makebox(0,0)[]{$D_1 $}}
\put(27.5,53){\makebox(0,0)[]{$A_2 $}}
\put(62.5,53){\makebox(0,0)[cc]{$B_2 $}}
\put(107.5,53){\makebox(0,0)[cc]{$C_2 $}}
\put(142.5,53){\makebox(0,0)[]{$D_2 $}}
\put(37.5,53){\makebox(0,0)[]{$A_{1,2} $}}
\put(72.5,53){\makebox(0,0)[cc]{$B_{1,2} $}}
\put(87.5,53){\makebox(0,0)[cc]{$C_{1,2} $}}
\put(122.5,53){\makebox(0,0)[]{$D_{1,2} $}}
\put(75,80){\circle*{2}}
\put(75,90){\circle*{2}}
\put(85,80){\circle*{2}}
\put(85,90){\circle*{2}}
\put(75,77){\makebox(0,0)[cc]{$\overline{d}$}}
\put(75,93){\makebox(0,0)[cc]{${d}$}}
\put(85,77){\makebox(0,0)[cc]{$\overline{p}$}}
\put(85,93){\makebox(0,0)[cc]{${p}$}}
\put(15,85){\makebox(0,0)[lc]{$A=A_1 \cup A_2 \cup A_{1,2} $}}
\put(15,75){\makebox(0,0)[lc]{$B=B_1 \cup B_2 \cup B_{1,2} $}}
\put(145,85){\makebox(0,0)[rc]{$C=C_1 \cup C_2 \cup C_{1,2} $}}
\put(145,75){\makebox(0,0)[rc]{$D=D_1 \cup D_2 \cup D_{1,2} $}}
\put(75,80){\line(0,1){10}}
\put(75,90){\line(1,-1){10}}
\put(85,80){\line(0,1){10}}
\put(85,90){\line(-1,-1){10}}
\end{picture}
}

\caption{Another construction for middle level sets ``$R$", which allows adjacent
ranks $k$, $k+1$ with equal rank $k$ and rank $k+1$ decks. Here, $\overline{d} $ is above all maximal elements of
$B\cup C$, {\em except} those that are circled.
$\overline{p} $ is above all maximal elements of
$B\cup C $, {\em except} those that are boxed.
By transitivity, $d$ and $p$ are above all elements of $B\cup C $. Moreover,
$d$
is above all maximal elements of
$A\cup D$, {\em except} those that are circled.
$p $ is above all maximal elements of
$A\cup D$, {\em except} those that are boxed.
}
\label{flatswit2}

\end{figure}

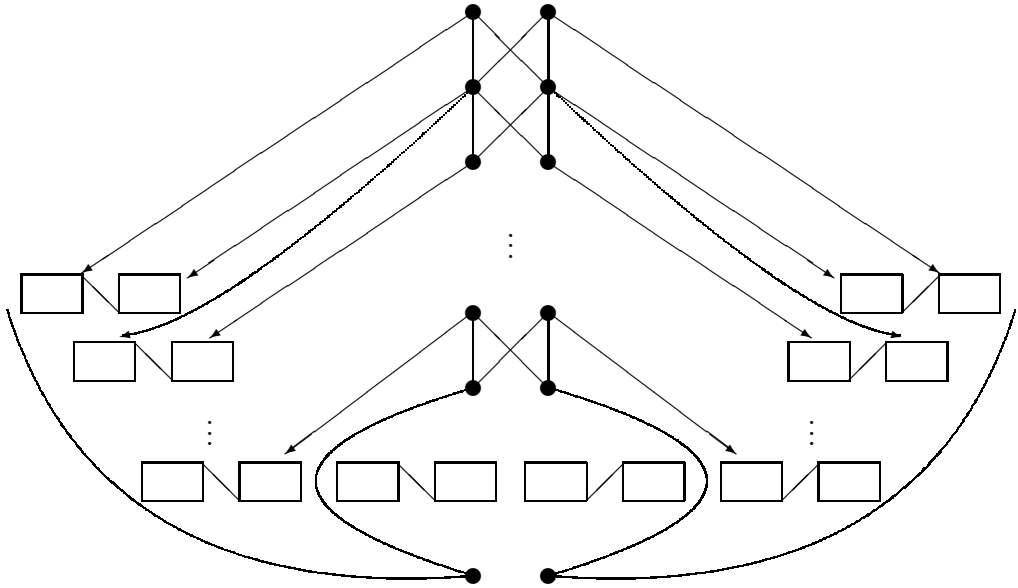
\begin{figure}


\unitlength 1.00mm 
\linethickness{0.4pt}
\ifx\plotpoint\undefined\newsavebox{\plotpoint}\fi 
\begin{picture}(137,81)(0,0)
\put(65,80){\line(0,-1){10}}
\put(65,70){\line(1,1){10}}
\put(75,80){\line(0,-1){10}}
\put(75,70){\line(-1,1){10}}
\put(65,80){\circle*{2}}
\put(75,80){\circle*{2}}
\put(65,70){\line(0,-1){10}}
\put(65,60){\line(1,1){10}}
\put(75,70){\line(0,-1){10}}
\put(75,60){\line(-1,1){10}}
\put(65,70){\circle*{2}}
\put(75,70){\circle*{2}}
\put(65,60){\circle*{2}}
\put(75,60){\circle*{2}}
\put(65,40){\line(0,-1){10}}
\put(65,30){\line(1,1){10}}
\put(75,40){\line(0,-1){10}}
\put(75,30){\line(-1,1){10}}
\put(65,40){\circle*{2}}
\put(75,40){\circle*{2}}
\put(65,30){\circle*{2}}
\put(75,30){\circle*{2}}
\put(70,50){\makebox(0,0)[cc]{$\vdots $}}
\put(47,15){\framebox(8,5)[cc]{}}
\put(55,20){\line(1,-1){5}}
\put(60,15){\framebox(8,5)[cc]{}}
\put(85,15){\framebox(8,5)[cc]{}}
\put(85,20){\line(-1,-1){5}}
\put(72,15){\framebox(8,5)[cc]{}}
\put(65,5){\circle*{2}}
\put(75,5){\circle*{2}}
\put(30,25){\makebox(0,0)[cc]{$\vdots $}}
\put(110,25){\makebox(0,0)[cc]{$\vdots $}}
\put(65,80){\vector(-3,-2){52}}
\put(65,70){\vector(-3,-2){38}}
\put(65,60){\vector(-3,-2){35}}
\put(65,40){\vector(-4,-3){25}}
\bezier{236}(65,70)(31.33,38)(18.33,37)
\put(19,37.1){\vector(-4,-1){1}}
\put(75,80){\vector(3,-2){52}}
\put(75,70){\vector(3,-2){38}}
\put(75,60){\vector(3,-2){35}}
\put(75,40){\vector(4,-3){25}}
\bezier{236}(75,70)(108.67,38)(121.67,37)
\put(121,37.1){\vector(4,-1){1}}
\bezier{352}(65,30)(23,17.67)(65,5)
\bezier{352}(75,30)(117,17.67)(75,5)
\bezier{364}(65,5)(15.33,1)(3,40.33)
\bezier{364}(75,5)(124.67,1)(137,40.33)
\put(21,15){\framebox(8,5)[cc]{}}
\put(29,20){\line(1,-1){5}}
\put(34,15){\framebox(8,5)[cc]{}}
\put(12,31){\framebox(8,5)[cc]{}}
\put(20,36){\line(1,-1){5}}
\put(25,31){\framebox(8,5)[cc]{}}
\put(5,40){\framebox(8,5)[cc]{}}
\put(13,45){\line(1,-1){5}}
\put(18,40){\framebox(8,5)[cc]{}}
\put(111,15){\framebox(8,5)[cc]{}}
\put(111,20){\line(-1,-1){5}}
\put(98,15){\framebox(8,5)[cc]{}}
\put(120,31){\framebox(8,5)[cc]{}}
\put(120,36){\line(-1,-1){5}}
\put(107,31){\framebox(8,5)[cc]{}}
\put(127,40){\framebox(8,5)[cc]{}}
\put(127,45){\line(-1,-1){5}}
\put(114,40){\framebox(8,5)[cc]{}}
\end{picture}

\caption{A construction to combine sets as in Figure \protect\ref{flatswit2}
to obtain sets in which many consecutive ranks have equal rank $k$ decks.
Fine structure of the involved sets is only indicated roughly
by showing with arrows which boxes (boxes are stand ins for the sets $A$, $B$, $C$, $D$)
go with which points. Crossovers ``point on left to box on right"
and symmetrically, as well as points $c_b $ and $c_t $ are omitted.
The bottom set (oval with some structure inside) can be a set $R$ or a set $Q$.
The two bottom points are lower bounds of all elements in the boxes.
None of the boxes have any comparabilities between them unless they are connected,
and the connection signifies comparability between the minimal and the maximal elements in
the indicated direction.
This implies that the elements of the four crown tower are the only elements
in this structure whose rank exceeds $2$.
}
\label{tower3}

\end{figure}

\begin{remark}
\label{consecstack}

{\rm
It may be considered unsatisfying that the ordered sets in
Lemma \ref{severalmiddle} are such that the ranks $k$ that yield equal rank $k$ decks are
separated from each other.
A small modification in the construction can also produce adjacent ranks
such that the rank $k$ decks of two nonisomorphic sets are equal.
The ordered set in Figure \ref{flatswit2}
has the same symmetry properties as the sets $R$
in Lemma \ref{severalmiddle}, except
that $\overline{d} $ and $\overline{p} $ are not minimal.

If these sets are now ``stacked" as indicated in Figure \ref{tower3}, with
$d$'s merged with $\overline{d} $'s and $p$'s merged with $\overline{p} $'s we obtain a
tower structure as in Lemma \ref{severalmiddle} without $Q$ and $\tilde{Q} $ attached.
The only difference is that the points that
eventually yield equal cards are concentrated in the four crown tower.
Replacing stretches of sets $R$ in a set as in Lemma \ref{severalmiddle} with
sets as in Figure \ref{tower3} introduces stretches in which
arbitrarily many consecutive ranks produce equal rank $k$ decks.
Also note that
the elements of the four crown tower in Figure \ref{tower3}
are the only elements with rank exceeding 2 and that we can again show that the
resulting sets will have equal rank $k$ neighborhood decks for all $k$.
}

\end{remark}

\section{A Folding Operation}

The construction that leads to Theorem \ref{tallones} produces
ordered sets that are ``tall"
in the sense that their height can exceed their width by an arbitrary factor.
Moreover, the ranks that produce equal rank $k$ decks are, with $2$ elements,
as small as can be.
It is possible to produce ``wider" examples by folding the examples in Theorem \ref{tallones}
appropriately.
This idea is explored in this section.

\begin{figure}

\centerline{
\unitlength 1.10mm 
\linethickness{0.4pt}
\ifx\plotpoint\undefined\newsavebox{\plotpoint}\fi 
\begin{picture}(130,41)(0,0)
\put(15,10){\line(1,-1){5}}
\put(20,5){\line(1,1){5}}
\put(25,10){\line(1,-1){5}}
\put(30,5){\line(1,1){5}}
\put(35,10){\line(0,1){10}}
\put(35,20){\line(-1,1){5}}
\put(30,25){\line(-1,-1){5}}
\put(25,20){\line(-1,1){5}}
\put(15,20){\line(0,-1){10}}
\put(30,5){\circle*{2}}
\put(20,5){\circle*{2}}
\put(32,25){\makebox(0,0)[lc]{\footnotesize $a_2 $}}
\put(18,25){\makebox(0,0)[rc]{\footnotesize $a_1 $}}
\put(15,30){\line(1,-1){5}}
\put(20,25){\line(1,1){5}}
\put(25,30){\line(1,-1){5}}
\put(30,25){\line(1,1){5}}
\put(35,35){\line(-1,1){5}}
\put(30,40){\line(-1,-1){5}}
\put(25,35){\line(-1,1){5}}
\put(20,40){\line(-1,-1){5}}
\put(30,25){\circle*{2}}
\put(20,25){\circle*{2}}
\put(30,40){\circle*{2}}
\put(20,40){\circle*{2}}
\put(20,25){\line(-1,-1){5}}
\put(25,32.67){\makebox(0,0)[cc]{$F$}}
\put(25,15){\makebox(0,0)[cc]{$S$}}
\put(10,40){\makebox(0,0)[cc]{$P$}}
\put(20,5){\line(-3,1){15}}
\put(30,5){\line(-4,1){20}}
\put(13,12){\makebox(0,0)[rc]{\tiny (further parts of $S$)}}
\put(30,40){\line(3,-1){15}}
\put(20,40){\line(4,-1){20}}
\put(37,33){\makebox(0,0)[lc]{\tiny (further parts of $F$)}}
\put(70,10){\line(1,-1){5}}
\put(75,5){\line(1,1){5}}
\put(80,10){\line(1,-1){5}}
\put(85,5){\line(1,1){5}}
\put(90,10){\line(0,1){10}}
\put(70,20){\line(0,-1){10}}
\put(85,5){\circle*{2}}
\put(75,5){\circle*{2}}
\put(80,15){\makebox(0,0)[cc]{$S$}}
\put(75,5){\line(-3,1){15}}
\put(85,5){\line(-4,1){20}}
\put(68,12){\makebox(0,0)[rc]{\tiny (further parts of $S$)}}
\put(100,27){\makebox(0,0)[cc]{\footnotesize $a_2 $}}
\put(90,27){\makebox(0,0)[cc]{\footnotesize $a_1 $}}
\put(100,25){\circle*{2}}
\put(90,25){\circle*{2}}
\put(100,25){\line(-2,-1){10}}
\put(70,20){\line(4,1){20}}
\put(90,25){\line(-2,-1){10}}
\put(80,20){\line(4,1){20}}
\put(120,10){\line(-1,-1){5}}
\put(115,5){\line(-1,1){5}}
\put(110,10){\line(-1,-1){5}}
\put(105,5){\line(-1,1){5}}
\put(105,5){\circle*{2}}
\put(115,5){\circle*{2}}
\put(110,13.33){\makebox(0,0)[cc]{$F$}}
\put(115,5){\line(3,1){15}}
\put(105,5){\line(4,1){20}}
\put(122,12){\makebox(0,0)[lc]{\tiny (further parts of $F$)}}
\put(110,11){\makebox(0,0)[cc]{\tiny (folded)}}
\put(80,12.33){\makebox(0,0)[cc]{\tiny (stayed)}}
\put(65,25){\makebox(0,0)[cc]{$P_{F,S} $}}
\put(20,30){\vector(0,1){5}}
\put(30,30){\vector(0,1){5}}
\put(20,10){\vector(0,1){10}}
\put(30,10){\vector(0,1){10}}
\put(75,10){\vector(0,1){10}}
\put(85,10){\vector(0,1){10}}
\put(115,15){\vector(0,-1){5}}
\put(105,15){\vector(0,-1){5}}
\put(15,35){\line(0,-1){5}}
\put(35,30){\line(0,1){5}}
\put(100,10){\line(0,1){5}}
\put(100,15){\line(-1,1){10}}
\put(90,25){\line(2,-1){20}}
\put(110,15){\line(-1,1){10}}
\put(100,25){\line(2,-1){20}}
\put(120,15){\line(0,-1){5}}
\end{picture}
}

\caption{The folding operation in Definition \protect\ref{foldoper}.}
\label{foldant}

\end{figure}
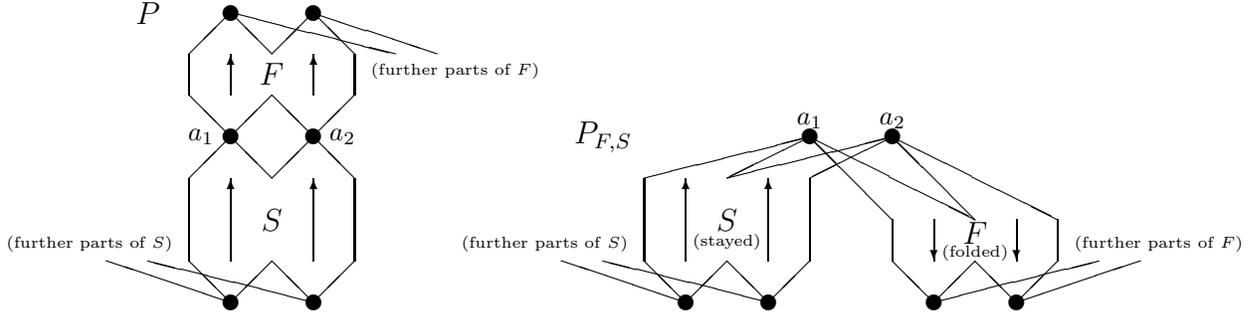

\begin{define}
\label{foldoper}

Let $P$ be a connected finite ordered set.
An antichain $A\subseteq P$ will be called a {\bf seam}
iff the removal of
\begin{enumerate}
\item
All comparabilities $x<y$
such that there are $a,b\in A$ with $x<a$ and $b<y$, and of
\item
The antichain $A$ and all comparabilities involving points of $A$,
\end{enumerate}
disconnects the ordered set $P$.
Call the resulting ordered set the {\bf $A$-separation of $P$}.
A seam will be called {\bf foldable}
iff
\begin{enumerate}
\item[F)]
For all $x,y\in P$ we have that if $x<y$ and there are $a,b\in A$ such that
$x<a$ and $b<y$, then there is a $c\in A$ such that $x<c<y$.
\end{enumerate}

A seam will be called {\bf breakable}
iff
\begin{enumerate}
\item[B)]
For all $x,y\in P$ we have that if there are $a,b\in A$ such that
$x<a$ and $b<y$, then $x<y$.
\end{enumerate}

Let $A$ be a foldable or breakable seam in
$P$, and let $F$
(the part to be folded)
and $S$ (the part that will stay as is) be nonempty unions of components
of the $A$-separation of $P$ such that
in $P$
no element of $S$ is above any element of $F$.
We define $P_{F,S} $
(also cf. Figure \ref{foldant})
to be the ordered set obtained from $P$ by
\begin{enumerate}
\item
Erasing all comparabilities $x<y$ in $P$
such that there are $a,b\in A$ with $x<a$ and $b<y$,

\item
Keeping the remaining comparabilities in $S\cup A$ as is,

\item
Reversing all comparabilities in $F\cup A$.

\end{enumerate}

(It is easy to see that this ``folding operation" produces an ordered set.)

\end{define}

The folding operation of Definition \ref{foldoper} is well behaved with respect to
isomorphism, isomorphism of cards and isomorphism of neighborhoods
as the next lemma shows.

\begin{lem}
\label{foldthem}

Let $P,P'$ be connected ordered sets and let $k\in {\mat N}$
be such that $A:=\{ a\in P: {\rm rank } _P (a)=k\} \subseteq P$ and
$A':=\{ a'\in P': {\rm rank } _{P'} (a')=k\} \subseteq P'$
are both foldable seams
(or both breakable seams)
in $P$ and $P'$ respectively.
Moreover assume that
we can choose
$F$, $S$, $F'$ and $S'$ as follows.

\begin{enumerate}
\item
$F$ and $S$ are nonempty unions of components
of the $A$-separation of $P$ such that
$F$ contains no points that are below any element of $A$ and
$S$ contains no points that are above any element of $A$.
\item
$F'$ and $S'$ are nonempty unions of components
of the $A'$-separation of $P'$ such that
$F'$ contains no points that are below any element of $A'$ and
$S'$ contains no points that are above any element of $A'$.
\item
\label{dualrank1}
All elements of $P$ that have dual rank $\geq k$ are in $S\cup A$
and all elements of $S$ have rank less than $k$.
\item
\label{dualrank2}
All elements of $P'$ that have dual rank $\geq k$ are in $S'\cup A'$
and all elements of $S$ have rank less than $k$.

\item
\label{nocrossover1}
No component of $F\cup A$ is isomorphic to the dual of a component of $S'\cup A'$.
\item
\label{nocrossover2}
No component of $F'\cup A'$ is isomorphic to the dual of a component of $S\cup A$.

\end{enumerate}

Then the following hold.

\begin{enumerate}
\item
\label{isoequiv}
$P$ is isomorphic to $P'$ iff $P_{F,S} $ is isomorphic to
$P' _{F',S'} $.
\item
\label{isofoldedcards}
For all points $p\in P$ and $p'\in P'$
such that
$P\setminus \{ p\} $ is isomorphic to
$P'\setminus \{ p'\} $
via an isomorphism that maps $A\setminus \{ p\} $ to $A'\setminus \{ p'\} $,
$S\setminus \{ p\} $ to $S'\setminus \{ p'\} $
and $F\setminus \{ p\} $ to $F'\setminus \{ p'\} $,
we have that
$P_{F,S} \setminus \{ p\} $ is isomorphic to
$P'_{F',S'} \setminus \{ p'\} $.
\item
\label{isofoldednhoods}
For all points $p\in P$ and $p'\in P'$ such that
${\rm rank}_P (p)={\rm rank}_{P'} (p')$
and
the neighborhood $\updownarrow _P p $
is isomorphic to the neighborhood
$\updownarrow _{P' } p'$, we have that
${\rm rank}_{P_{F,S}} (p)={\rm rank}_{P'_{F',S'}} (p')$
and the neighborhood
$\updownarrow _{P_{F,S} } p $
is isomorphic to the neighborhood
$\updownarrow _{P'_{F',S'} } p'$.
\end{enumerate}

\end{lem}

{\bf Proof.}
For part \ref{isoequiv} first let $P_{F,S} $ be isomorphic to $P'_{F',S'} $
via the isomorphism $\Phi :P_{F,S} \to P'_{F',S'} $.
By hypotheses \ref{dualrank1} and \ref{dualrank2}
we have that $A=\{ p\in P_{F,S} : {\rm rank } _{P_{F,S} } (p)=k\} $ and
$A'=\{ p'\in P'_{F,S} : {\rm rank } _{P'_{F,S} } (p')=k\} $.
Therefore $\Phi [A]=A'$.
Then by hypotheses
\ref{nocrossover1}
and
\ref{nocrossover2}
we have that $\Phi |_{S\cup A} $ is an isomorphism between $S\cup A$ and $S'\cup A'$
and $\Phi |_{F\cup A} $ is an isomorphism between $F\cup A$ and $F'\cup A'$, where the
respective sets carry the orders that are induced by
$P_{F,S} $ and by $P' _{F',S'} $, respectively.
Since the folding construction only reverses the orders on $F\cup A$ and $F'\cup A'$,
$\Phi |_{S\cup A} $ is an isomorphism between $S\cup A$ and $S'\cup A'$
and $\Phi |_{F\cup A} $ is an isomorphism between $F\cup A$ and $F'\cup A'$, where the
respective sets carry the orders that are induced by $P$ and by
$P'$, respectively.
Now if both $A$ and $A'$ are foldable seams, then all comparabilities
between elements of $F$ and $S$ (and of $F'$ and $S'$) are induced by transitivity through an
element of $A$ ($A'$ respectively).
This implies that $\Phi $ is an isomorphism between $P$ and $P'$.
If both $A$ and $A'$ are breakable seams, then all comparabilities
between elements of $F$ and $S$ (and of $F'$ and $S'$) that are related to any element of
$A$ ($A'$ respectively) are present in $P$ ($P'$ respectively) and these are all comparabilities
between elements of $F$ and $S$ (and of $F'$ and $S'$).
Again $\Phi $ is an isomorphism between $P$ and $P'$.

For the converse, let $P $ be isomorphic to $P' $
via the isomorphism $\Psi :P \to P' $.
Then
$\Psi |_{S\cup A} $ is an isomorphism between $S\cup A$ and $S'\cup A'$
and $\Psi |_{F\cup A} $ is an isomorphism between $F\cup A$ and $F'\cup A'$, where the
respective sets carry the orders that are induced by $P$ and by
$P'$, respectively.
Since there are no comparabilities between elements of $F$ and $S$ ($F'$ and $S'$ respectively)
in $P_{F,S} $ ($P'_{F',S'} $), this means that
$\Psi $ is also an isomorphism between $P_{F,S} $ and $P'_{F',S'} $.

For part \ref{isofoldedcards},
let $\Psi :P\setminus \{ p\} \to P\setminus \{ p'\} $ be an isomorphism
that maps $A\setminus \{ p\} $ to $A'\setminus \{ p'\} $, $S\setminus \{ p\} $ to $S'\setminus \{ p'\} $
and $F\setminus \{ p\} $ to $F'\setminus \{ p'\} $.
Then $\Psi $ is a bijection between
$P_{F,S} \setminus \{ p\} $ and $P'_{F',S'} \setminus \{ p'\} $
that is an isomorphism between
$(F\cup A)\setminus \{ p\} $ and
$(F'\cup A')\setminus \{ p'\} $ and between
$(S\cup A)\setminus \{ p\} $ and
$(S'\cup A')\setminus \{ p'\} $, respectively.
Thus
$\Psi $ is an isomorphism between
$P_{F,S} \setminus \{ p\} $ and
$P'_{F',S'} \setminus \{ p'\} $.

For part \ref{isofoldednhoods} first note that
because the neighborhoods are isomorphic,
the dual rank of
$p$ in $P$ is equal to the dual rank of $p'$ in $P'$.
Since $p$ and $p'$ have the same rank in their respective unfolded sets,
their ranks in the folded sets is either their
original rank or their original dual rank.
Either way, their ranks in the folded sets are equal.
If the original rank of $p$ and $p'$ is less than $k$, then
their neighborhoods in the folded sets are obtained by
discarding all elements of rank greater than $k$ from the
original neighborhoods.
If the original rank of $p$ and $p'$ is greater than $k$, then
their neighborhoods in the folded sets are obtained by
discarding all elements of rank less than $k$ from the
original neighborhoods and dualizing the order.
If $p\in A$ and $p'\in A'$, then their neighborhoods
in the folded sets are obtained
from the original neighborhoods by folding at
$p$ or $p'$, respectively.
In all cases the neighborhoods in the folded sets are isomorphic.
\qed

\begin{theorem}
\label{populatedlevels}

For every sequence of natural numbers $\{ s_n \} _{n\geq 1} $,
there is a sequence of pairs of ordered sets $(Q_1 ^n , Q_2 ^n )$
such that
\begin{enumerate}
\item
$Q_1 ^n $ is not isomorphic to $Q_2 ^n $,
\item
\label{eqmaxmindeck2}
$Q_1 ^n $ and $Q_2 ^n $ have equal maximal and minimal decks,
\item
\label{eqrankkdecks2}
There are ranks $r_0, \ldots , r_{n} $ such that $Q_1 ^n $ and $Q_2 ^n $ have
equal marked rank $r_i $ decks,
and the ranks $r_1 , \ldots ,r_n $ do not contain
extremal elements,
\item
\label{eqcardeachrnk}
For each rank $k$ there is at least one element of rank $k$ in $Q_1 ^n $ and $Q_2 ^n $
such that
removal of these elements produces isomorphic marked rank $k$ cards,

\item
\label{eqcardeqremnhood}
For all the above mentioned isomorphic cards $Q_1 ^n \setminus \{ x\} $
and $Q_2 ^n \setminus \{ x\} $ the
neighborhoods $\updownarrow _{Q_1 ^n } x $ and $\updownarrow _{Q_2 ^n } x $
are isomorphic.

\item
$\displaystyle{ \liminf _{n\to \infty } ECR(Q_1 ^n , Q_2 ^n )\geq {1\over 10}} $,

\item
\label{equalnhooddecksforQ}
For every rank $k$, the rank $k$ neighborhood decks of $Q_1 ^n $ and $Q_2 ^n $
are equal.

\item
$Q_1 ^n $ and $Q_2 ^n $ have at least $s_n $ maximal elements, at least
$s_n $ minimal elements and for each $k_i $ at least $s_n $ elements of rank $k_i $.
\end{enumerate}

\end{theorem}

{\bf Proof.} This result is a direct consequence of Theorem \ref{tallones} and
Lemma \ref{foldthem}.
Without loss of generality, assume that $s_n $ is even
(otherwise replace it with $s_n +1$).
Now let $\displaystyle{ t_n :={s_n \over 2}(n+2) } $,
consider the pair of sets
$\displaystyle{\left( H^{s_n \over 2} , K^{s_n \over 2} \right) :=\left( P_1 ^{t_n } , P_2 ^{t_n } \right) } $
from Theorem \ref{tallones},
and let $k_0 , \ldots , k_{t_n } $ be the ranks mentioned in Theorem \ref{tallones}.
Now
for $\displaystyle{ i \in \left\{ {s_n \over 2}, \ldots , 1\right\} } $,
obtain $\left( H^{i-1} , K^{i-1} \right) $ from $\left( H^{i} , K^{i} \right) $ by folding the sets
$H^{i} , K^{i} $ as indicated in Lemma \ref{foldthem} at rank
$[{\rm height} (Q)+i(n+2) {\rm height}(R)]$.
At each fold, the
height of the part $S$ that stays as is is greater than the height of the part $F$ that is folded and
both parts are connected. Thus there cannot be any isomorphisms between duals and
Lemma \ref{foldthem} applies without a problem.

In a last step, apply the dual of Lemma \ref{foldthem} at rank $[{\rm height} Q]$
(which is possible because the height of $Q$ can be chosen to be equal to the height of $R$).

Let $(Q_1 ^n , Q_2 ^n )$ be the pair of ordered sets thus obtained.
By Lemma \ref{foldthem}, the pair $(Q_1 ^n , Q_2 ^n )$ is as desired.
The last step of folding up the bottom set $Q$ guarantees
part \ref{eqcardeachrnk}.
\qed

\begin{remark}

{\rm
While the examples in Theorem \ref{populatedlevels} are not counterexamples to the reconstruction conjecture,
they show that the subtlety of order reconstruction reaches beyond tools available today.
All parameters that we know to be reconstructible are equal for these sets.
Moreover, several parameters, such as maximal and minimal decks as well as rank $k$ decks, plus appropriate
markings, are equal also.
These parameters
have not yet been proven to be reconstructible.
}

\end{remark}

\begin{remark}
\label{allbut2}

{\rm
Using sets as in Remark \ref{consecstack} and their duals and
folding appropriately, using Lemma \ref{foldthem} and its dual,
it is now possible to produce pairs of nonisomorphic ordered sets
of arbitrary height
for
which all ranks except ranks $1$ and $2$ produce equal rank $k$ decks, for which all
ranks have arbitrarily many elements and for which even in ranks $1$ and $2$ there
are arbitrarily many elements whose removal produces isomorphic cards.
The equal card ratio will approach at least $\displaystyle{ 1\over 10} $
with the constructions
available in this paper.
}

\end{remark}

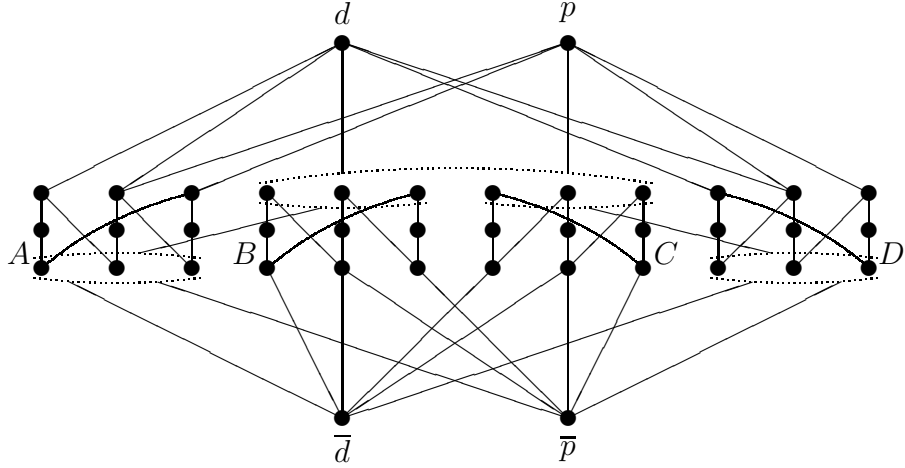
\begin{figure}

\centerline{
\unitlength 1.00mm
\linethickness{0.4pt}
\begin{picture}(123.00,64.00)
\put(10.00,30.00){\circle*{2.00}}
\put(20.00,30.00){\circle*{2.00}}
\put(30.00,30.00){\circle*{2.00}}
\put(40.00,30.00){\circle*{2.00}}
\put(50.00,30.00){\circle*{2.00}}
\put(60.00,30.00){\circle*{2.00}}
\put(70.00,30.00){\circle*{2.00}}
\put(80.00,30.00){\circle*{2.00}}
\put(90.00,30.00){\circle*{2.00}}
\put(100.00,30.00){\circle*{2.00}}
\put(110.00,30.00){\circle*{2.00}}
\put(120.00,30.00){\circle*{2.00}}
\put(10.00,40.00){\circle*{2.00}}
\put(20.00,40.00){\circle*{2.00}}
\put(30.00,40.00){\circle*{2.00}}
\put(40.00,40.00){\circle*{2.00}}
\put(50.00,40.00){\circle*{2.00}}
\put(60.00,40.00){\circle*{2.00}}
\put(70.00,40.00){\circle*{2.00}}
\put(80.00,40.00){\circle*{2.00}}
\put(90.00,40.00){\circle*{2.00}}
\put(100.00,40.00){\circle*{2.00}}
\put(110.00,40.00){\circle*{2.00}}
\put(120.00,40.00){\circle*{2.00}}
\put(10.00,40.00){\line(1,-1){10.00}}
\put(20.00,30.00){\line(0,1){10.00}}
\put(20.00,40.00){\line(1,-1){10.00}}
\put(30.00,30.00){\line(0,1){10.00}}
\put(10.00,30.00){\line(0,1){10.00}}
\put(120.00,40.00){\line(-1,-1){10.00}}
\put(110.00,30.00){\line(0,1){10.00}}
\put(110.00,40.00){\line(-1,-1){10.00}}
\put(100.00,30.00){\line(0,1){10.00}}
\put(120.00,30.00){\line(0,1){10.00}}
\put(50.00,60.00){\circle*{2.00}}
\put(80.00,60.00){\circle*{2.00}}
\put(10.00,40.00){\line(2,1){40.00}}
\put(50.00,60.00){\line(-3,-2){30.00}}
\put(20.00,40.00){\line(3,1){60.00}}
\put(80.00,60.00){\line(-5,-2){50.00}}
\put(120.00,40.00){\line(-2,1){40.00}}
\put(80.00,60.00){\line(3,-2){30.00}}
\put(110.00,40.00){\line(-3,1){60.00}}
\put(50.00,60.00){\line(5,-2){50.00}}
\put(40.00,40.00){\line(1,-1){10.00}}
\put(50.00,30.00){\line(0,1){10.00}}
\put(50.00,40.00){\line(1,-1){10.00}}
\put(60.00,30.00){\line(0,1){10.00}}
\put(40.00,30.00){\line(0,1){10.00}}
\put(90.00,40.00){\line(-1,-1){10.00}}
\put(80.00,30.00){\line(0,1){10.00}}
\put(80.00,40.00){\line(-1,-1){10.00}}
\put(70.00,30.00){\line(0,1){10.00}}
\put(90.00,30.00){\line(0,1){10.00}}
\put(50.00,10.00){\circle*{2.00}}
\put(80.00,10.00){\circle*{2.00}}
\put(40.00,30.00){\line(1,-2){10.00}}
\put(50.00,10.00){\line(0,1){20.00}}
\put(50.00,30.00){\line(3,-2){30.00}}
\put(80.00,10.00){\line(-1,1){20.00}}
\put(90.00,30.00){\line(-1,-2){10.00}}
\put(80.00,10.00){\line(0,1){20.00}}
\put(80.00,30.00){\line(-3,-2){30.00}}
\put(50.00,10.00){\line(1,1){20.00}}
\bezier{70}(39.00,41.30)(65.00,45.30)(91.00,41.30)
\put(50.00,60.00){\line(0,-1){17.37}}
\put(80.00,60.00){\line(0,-1){17.37}}
\bezier{30}(39.00,38.70)(50.00,37.30)(61.00,38.70)
\bezier{30}(9.00,31.30)(20.00,32.70)(31.00,31.30)
\put(50.00,10.00){\line(-2,1){36.33}}
\put(80.00,10.00){\line(-3,1){54.50}}
\put(47.00,38.00){\line(-4,-1){24.00}}
\bezier{30}(9.00,28.70)(20.00,27.30)(31.00,28.70)
\bezier{30}(91.00,38.70)(80.00,37.30)(69.00,38.70)
\bezier{30}(121.00,31.30)(110.00,32.70)(99.00,31.30)
\put(80.00,10.00){\line(2,1){36.33}}
\put(50.00,10.00){\line(3,1){54.50}}
\put(83.00,38.00){\line(4,-1){24.00}}
\bezier{30}(121.00,28.70)(110.00,27.30)(99.00,28.70)
\put(10.00,35.00){\circle*{2.00}}
\put(20.00,35.00){\circle*{2.00}}
\put(30.00,35.00){\circle*{2.00}}
\put(40.00,35.00){\circle*{2.00}}
\put(50.00,35.00){\circle*{2.00}}
\put(60.00,35.00){\circle*{2.00}}
\put(70.00,35.00){\circle*{2.00}}
\put(80.00,35.00){\circle*{2.00}}
\put(90.00,35.00){\circle*{2.00}}
\put(100.00,35.00){\circle*{2.00}}
\put(110.00,35.00){\circle*{2.00}}
\put(120.00,35.00){\circle*{2.00}}
\bezier{192}(10.00,30.00)(18.00,37.00)(30.00,40.00)
\bezier{192}(40.00,30.00)(48.00,37.00)(60.00,40.00)
\bezier{192}(120.00,30.00)(112.00,37.00)(100.00,40.00)
\bezier{192}(90.00,30.00)(82.00,37.00)(70.00,40.00)
\put(50.00,6.00){\makebox(0,0)[cc]{$\overline{d} $}}
\put(80.00,6.00){\makebox(0,0)[cc]{$\overline{p} $}}
\put(50.00,64.00){\makebox(0,0)[cc]{${d} $}}
\put(80.00,64.00){\makebox(0,0)[cc]{${p} $}}
\put(7.00,32.00){\makebox(0,0)[cc]{$A$}}
\put(37.00,32.00){\makebox(0,0)[cc]{$B$}}
\put(93.00,32.00){\makebox(0,0)[cc]{$C$}}
\put(123.00,32.00){\makebox(0,0)[cc]{$D$}}
\end{picture}
}

\caption{An ordered set $R$ as needed in Lemma \protect\ref{severalmiddle}.
Connected dotted arches again indicate a complete bipartite structure.
The sets in this figure can be used to construct sets with equal maximal and minimal
decks for which the ratio of the number of extremal elements to the size of the set is as large as
possible to date.
}
\label{switch4}

\end{figure}

\begin{remark}

{\rm
Aside from results in this paper, the only construction to obtain sets with
equal maximal decks with $k>3$ cards is due to
Ille and Rampon (cf. \cite{JXsurvey}, Section 8.2.2).
The size of their sets is exponential in the number of maximal elements.
The size of the sets in
Theorem \ref{populatedlevels} is linear in the number of maximal
elements.
In this construction, to date the largest ratio of extremal elements to the size of the set is
achieved with sets $R$ as in Figure \ref{switch4}. The ratio approaches $\displaystyle{1\over 19} $ and it is
achieved by folding at every merge of sets $R$ and at the merge of the top set $R$ with $\tilde{Q} $.
With sets as presented earlier, the ratio approaches $\displaystyle{ 1\over 20} $.
}

\end{remark}

\section{Further Consequences of the Folding Construction}

Aside from giving access to examples on reconstruction, the folding construction of Definition \ref{foldoper}
also allows access to some results that can simplify the start of a reconstruction proof.

Recall that two elements $x,y$ of an ordered set are called {\bf adjacent} iff
$x<y$ and for all $z\in P$ with $x\leq z\leq y$ we have $z\in \{ x,y\} $
(or the dual of this statement).
An ordered set $P$ is called {\bf graded} iff there is a function
$g:P\to {\mat N} $ such that if $x,y\in P$ are adjacent then $|g(x)-g(y)|=1$.
If the function $g$ can be chosen to be the function that assigns each element its rank, we say $P$
is {\bf graded by the rank function}.

\begin{prop}
\label{shrinkheight}

If ordered sets of height $2$ are reconstructible, then
all ordered sets that are graded by the rank function are reconstructible from their
marked ranked decks.

\end{prop}

{\bf Proof.}
Suppose $P $ and $P' $ were two nonisomorphic ordered sets that are graded by their
rank functions and which have equal marked ranked decks.
Note that in an ordered set that is graded by the rank function, every rank $P_k $ is
a foldable seam.
Repeatedly folding $P $ and $P' $ as indicated in
Lemma \ref{foldthem}, at the rank that is one less than the height of the set
and stopping when the
resulting sets have height $2$, would create
two nonisomorphic ordered sets of height $2$ with the same deck.
(Equality of the marked ranked decks is needed here to make sure that all
hypotheses in part \ref{isofoldedcards} of the conclusion of Lemma \ref{foldthem}
are satisfied.)
We have a contradiction to the assumption that
ordered sets of height $2$ are reconstructible.
\qed

\begin{remark}
\label{ht1remark}

{\rm
Note that except for extremely symmetric sets, we could actually fold
until we have an ordered set of height $1$ in the proof of Proposition \ref{shrinkheight}.
This observation underscores once more the importance of understanding more about the
reconstruction of ordered sets of height $1$.
}

\end{remark}

\begin{prop}
\label{onlyoneseam}

If all ordered sets with at most one foldable or breakable
seam are reconstructible from their
marked ranked decks, then
all ordered sets are reconstructible from their
marked ranked decks

\end{prop}

{\bf Proof.}
Suppose $P_1 $ and $P_2 $ are two nonisomorphic ordered sets that have equal marked ranked decks.
Find the two consecutive seams for which the difference in the ranks is maximal.
Then fold all seams except these two.
In a final step, one of these two can be folded also.
The result are two nonisomorphic ordered sets with equal marked ranked decks and
at most one foldable or breakable seam.
\qed

\begin{remark}

{\rm Proposition \ref{onlyoneseam} is reminiscent of the result in \cite{Yongzhi}.
Removal of a foldable seam disconnects the
covering graph, while removal of a breakable seam is a step towards
disconnecting it (cross-connections still need to be erased).
This means for reconstruction work, some type of ``vertical connectivity"
could be assumed in any attempt to reconstruct ordered sets.
Unlike the result in \cite{Yongzhi}, there is no restriction on the size of the separating set.
}

\end{remark}

\begin{remark}

{\rm
Similar to Remark \ref{ht1remark}, in the proof of Proposition \ref{onlyoneseam},
unless we encounter a highly symmetrical situation, both seams
can be folded at the end of the proof, leaving ordered sets without seams.
}

\end{remark}

\section{The Role of Rigidity}
\label{rigsect}

We have seen that the rigid bottom and top sets $Q$ and $\tilde{Q} $ play a crucial role in the
development of the examples presented here.
Without these rigid ``anchors" the sets $P_1 $ and $P_2 $ would ``untwist" and be isomorphic to each other.
Thus it is reasonable to shed some light on the role of rigidity in reconstruction.
The following results show that
isomorphism between various types of rigid substructures immediately leads to an isomorphism between the
sets.
Such rigid substructures are often recognizable from a maximal card or would be
recognizable if the marked ranked deck is available.
Thus the presence of such structures leads to reconstructibility.
Consequently, if there is a counterexample to the order reconstruction conjecture,
then it must be made up largely of non-rigid structures.
For other properties that a counterexample must have, cf. \cite{Schrigid,Schmoreexam}.

In particular, the results in this section seem to indicate that any
development of counterexamples based on the present examples will
need to

\begin{itemize}
\item
Remove the rigid ``anchors" without letting the two sets become isomorphic,
\item
Avoid the introduction of rigid structures in middle levels.
\end{itemize}

To state the results in this section, for an ordered set $P$ we denote

\begin{eqnarray*}
P_k & := & \{ x\in P: {\rm rank} (x)= k\} \\
P_{k\downarrow } & := & \{ x\in P: {\rm rank} (x)\leq k\} \\
P_{k\uparrow } & := & \{ x\in P: {\rm rank} (x)\geq k\} .
\end{eqnarray*}

%
%
%

\begin{prop}
\label{rigidtop}

Let $k\geq 0$ and let $P,Q$ be ordered sets with equal marked rank $k$ decks.
If
\begin{enumerate}
\item
$|P_k |>1$ and $|Q_k |>1 $, and
\item
No two elements of rank $k$
have the same strict upper and lower bounds
in either $P$ or $Q$, and
\item
$P\setminus P_k $ and $Q\setminus Q_k $ are both rigid,
\end{enumerate}
then $P$ is isomorphic to $Q$.

\end{prop}

{\bf Proof.}
Let $P\setminus \{ p\} $ and $Q\setminus \{ q\} $ be isomorphic rank $k$ cards and
let $\psi :P\setminus \{ p\} \to Q\setminus \{ q\} $ be an isomorphism that
preserves the original rank of each element.
Let $P\setminus \{ p'\} $ and $Q\setminus \{ q'\} $ be isomorphic rank $k$ cards
with $p\not= p'$, $q\not= q'$ and
let $\phi :P\setminus \{ p'\} \to Q\setminus \{ q'\} $ be an isomorphism that
preserves the original rank of each element.
Then $\phi |_{P\setminus P_k } =\psi |_{P\setminus P_k } $.

If $\psi (x)=q'$ and $x\not= p'$, then $\psi (p')\not= q'$, which means $y:=\phi ^{-1} (\psi (p'))\not= p'$.
This means that $y$ and $p'$ have the same sets of strict upper and lower bounds, a contradiction.
Thus $\psi (p')=q'$ and symmetrically $\phi (p)=q$.

Now for any element $x\in P_k \setminus \{ p,p'\} $, the inequality $\psi (x)\not= \phi (x)$ would imply
$\phi ^{-1} (\psi (x))\not= x$ has the same strict upper and lower bounds as $x$, a contradiction.
Thus
for all $x\in P_k \setminus \{ p,p'\} $, we have $\psi (x)= \phi (x)$.

Via $\phi |_{P\setminus P_k } =\psi |_{P\setminus P_k } $ we immediately conclude that
$$\Phi (x):=\cases{
\psi (x); & for $x\not= p$, \cr
\phi (p); & for $x=p$, \cr
} $$
is an isomorphism between $P$ and $Q$.
\qed

\begin{define}

Let $P$ be an ordered set and let $0<k<l$.
The
subset $P_{k\uparrow } \cap P_{l\downarrow } $ of $P$
is called a {\bf rigid separator} iff
$P_{k\uparrow } \cap P_{l\downarrow } $
is rigid, there are elements in $P$ of rank $>l$ and
no element of rank $<k$ is a lower cover of an element of rank $>l$.

\end{define}

\begin{prop}
\label{rigsep}

Let $P$ and $Q$ be ordered sets and let $0<k<l$. If
\begin{enumerate}
\item
$P_{k\uparrow } \cap P_{l\downarrow } $
and $Q_{k\uparrow } \cap Q_{l\downarrow } $ are rigid separators,
\item
There is an isomorphism $\phi :P\setminus \{ M_P \} \to Q\setminus \{ M_Q \}  $,
where
$M_P $ and $M_Q $ denote maximal elements of rank $r>l$ in $P$ and $Q$, respectively, and where
${\rm rank} _Q (\phi (x)) ={\rm rank} _P (x)$ for all $x$,
\item
There is an isomorphism
$\psi :P\setminus \{ m_P \} \to Q\setminus \{ m_Q \} $, where
$m_P $ and $m_Q $ denote
minimal elements in $P$ and $Q$, respectively,
and where
${\rm rank} _Q (\psi (x)) ={\rm rank} _P (x)$ for all $x$,
\end{enumerate}
then
$P$ is isomorphic to $Q$.

\end{prop}

{\bf Proof.}
Because we must have
$\phi |_{P_{k\uparrow } \cap P_{l\downarrow } } =\psi |_{P_{k\uparrow } \cap P_{l\downarrow } } $
and because there are no adjacencies that ``cross" $P_{k\uparrow } \cap P_{l\downarrow } $
the function
$$\Phi (x):=\cases{
\psi (x); & if ${\rm rank} (x)>l$, \cr
\phi (x); & if ${\rm rank} (x)\leq l$, \cr
} $$
is an isomorphism between $P$ and $Q$.
\qed

\section{Conclusion}

The examples presented in this paper show that even substantial partial information
on the deck of an ordered set is not sufficient to
effect reconstruction. In particular (see Theorem \ref{populatedlevels}),
the maximal deck plus the minimal deck plus
$n+1$ rank $k$ decks, plus the rank $k$ neighborhood decks
are not sufficient for reconstruction
even if all these decks have substantially more than $2$ cards.
Moreover, by Remark \ref{allbut2} even all rank $k$ decks except for the ones for
rank $1$ and rank $2$ are not sufficient to effect reconstruction. Again the
absolute number of cards within these ranks is immaterial.

Future reconstruction research has to take these facts into account.
Proof attempts that do not consider enough information from the deck will not succeed.
On the other hand, the examples presented here show what kind of
information might effect reconstruction or lead to a counterexample.

\begin{enumerate}
\item
All examples of ordered sets with equal maximal, minimal and
some equal rank $k$ decks so far have ``small waists".
That is, the ranks that have equal decks all have fewer elements
(by at least a factor $6$, even if we use sets $R$ as in Figure \ref{switch4})
than other ranks nearby.
It thus should be instructive to investigate what can be
concluded from equality of rank $k$ decks, where $k$ is such that
no other rank has more elements than the $k^{\rm th} $ rank.

\item
Alternatively, the construction of examples could be expanded to the
point where the overlap between the decks of two nonisomorphic sets
becomes so large that the decks would indeed have to be equal.
The present examples have been shown to be quite malleable. They might point the
way towards examples with similarly strong properties and
larger equal card ratios or maybe even a counterexample
overall. The largest equal card ratios in ordered sets observed so far slightly exceed 50\%
(see \cite{JHthes}), but the examples do not have the equal subdecks that the examples presented here have.

\item
Along these lines, would it be true that
if there is a sequence of pairs of nonisomorphic ordered sets $(P_1 ^n , P_2 ^n )$
such that $ECR (P_1 ^n , P_2 ^n ) \to 1$ as $n\to \infty $, then
there is a counterexample to the reconstruction conjecture?

\item
On p. 185 of \cite{Manv}, P. Stockmeyer is quoted to have said
only {\em half} in jest that ``The reconstruction conjecture [for graphs]
is not true, but the smallest counterexample has 87 vertices and will
never be found." The present examples seem to show that if there is a
counterexample to the order reconstruction conjecture, it would have to be
of substantial size and complexity.
At the same time, the approach of analyzing macrostructure
as in Lemma \ref{severalmiddle}
and microstructure as in Lemma \ref{getsmallerR}
separately may allow to
break down the complexity to manageable stages.

\item
Finally, Proposition \ref{onlyoneseam} shows that we can concentrate on
ordered sets that have a certain type of ``vertical connectivity".

\end{enumerate}

\end{document}